\documentclass{article}
\usepackage{amssymb}
\usepackage{amsmath}
\usepackage{amsfonts}

\usepackage{graphicx}
\usepackage[bookmarks=false]{hyperref}
\begin{document}

\author{Jean Petitot}
\title{Singularity theory and heat equation. \\
 \bigskip
 Th\'{e}orie des singularit\'{e}s et \'{e}quations de diffusion: 
illustration graphique d'un exemple.}
\date{21 avril 2004}
\maketitle

\begin{abstract}
The paper (in French) exemplifies graphically a solution of the heat equation which is a 1-dimensional unfolding of an elliptic umbilic catastrophe. The example is due to James Damon and adapts Thom-Mather's singularity theory to multiscale models of scale-space analysis in image processing.

\bigskip

Nous illustrons graphiquement un exemple, d\^{u} \`{a} James Damon, de
solution de l'\'{e}quation de la chaleur qui est un d\'{e}ploiement
\`{a} un param\`{e}tre d'un ombilic elliptique. Il fait partie d'une
adaptation de la th\'{e}orie des singularit\'{e}s de Thom-Mather aux
mod\`{e}les multi\'{e}chelle de l'analyse d'images (scale-space analysis).
\end{abstract}

\section{Introduction}

Pour \'{e}valuer l'h\'{e}ritage de l'\oe uvre de Ren\'{e} Thom, il faut je
crois voir les choses \`{a} un double niveau. Il y a d'abord la
post\'{e}rit\'{e} directe de ses id\'{e}es et puis il y a ensuite le
d\'{e}veloppement, dans d'autres contextes, d'id\'{e}es de type thomien
concernant l'application \`{a} des domaines originaux d'un certain nombre de
formalismes de g\'{e}om\'{e}trie diff\'{e}rentielle, de th\'{e}orie des
singularit\'{e}s, de th\'{e}orie des bifurcations, de ruptures de
sym\'{e}trie, de ph\'{e}nom\`{e}nes critiques, de ph\'{e}nom\`{e}nes de
transitions de phases qui, jusque-l\`{a}, ne concernaient que la physique.
Ren\'{e} Thom a \'{e}t\'{e} l'un des principaux fondateurs d'un nouveau
paradigme, le paradigme que j'appelle ``morphodynamique'' (th\'{e}orie
dynamique des formes), pour la biologie th\'{e}orique et les sciences
cognitives et ce dernier est devenu dominant dans de nombreux domaines.

En ce qui concerne la post\'{e}rit\'{e} directe des id\'{e}es de Ren\'{e}
Thom, il y a eu bien s\^{u}r le travail de ses disciples. M\^{e}me si
Ren\'{e} Thom n'a pas constitu\'{e} d'\'{e}cole \`{a} proprement parler, il
a f\'{e}cond\'{e} la pens\'{e}e d'un nombre consid\'{e}rable de personnes,
aussi bien dans le domaine des math\'{e}matiques pures que dans celui des
applications et celui de la philosophie. Tous les chercheurs qu'il a
inspir\'{e}s ont continu\'{e} \`{a} d\'{e}velopper certains aspects de son 
\oe uvre m\^{e}me si, sur le plan sociologique et m\'{e}diatique, apr\`{e}s
les grands d\'{e}bats des ann\'{e}es 70, la th\'{e}orie des catastrophes con%
\c{c}ue au sens strict a connu un certain retrait. Mais ce retrait,
d'ailleurs tout relatif, a \'{e}t\'{e}, selon moi, tr\`{e}s largement
compens\'{e} par l'extraordinaire succ\`{e}s d'id\'{e}es \textit{de type}
thomien. En particulier dans le domaine des sciences cognitives que je
connais un peu, les id\'{e}es introduites par Thom et Zeeman \`{a} la fin
des ann\'{e}es 60 sont devenues des id\'{e}es-force.

L'un des centres d'int\'{e}r\^{e}t majeurs de Ren\'{e} Thom a \'{e}t\'{e}
d'\'{e}laborer une th\'{e}orie des formes qui puisse \^{e}tre compatible
\`{a} la fois avec la gen\`{e}se objective de morphologies \`{a}
l'int\'{e}rieur des substrats mat\'{e}riels physico-chimiques et biologiques
(cf. ses travaux fondamentaux sur l'embryogen\`{e}se) et aussi avec la
perception subjective de morphologies \`{a} partir du traitement cognitif du
flux optique manifestant ces substrats.\footnote{%
Le lecteur interess\'{e} par les diff\'{e}rents aspects de la th\'{e}orie de
la forme chez Ren\'{e} Thom pourra consulter notre article ``Forme'' dans l'%
\textit{Encyclop\ae dia Universalis }(Petitot [1989]).}

L'un des premiers principes de la m\'{e}taphysique thomienne de la nature
est que les \textit{singularit\'{e}s} constituent, en tant que
discontinuit\'{e}s qualitatives, \textit{l'interface ph\'{e}nom\'{e}nologique%
} entre le monde objectif et la conscience perceptive. Cela \'{e}tait
d\'{e}j\`{a} bien connu pour l'espace depuis l'esth\'{e}tique
transcendantale kantienne, l'espace op\'{e}rant \`{a} la fois comme cadre de
l'objectivit\'{e} physique et comme comme forme de la manifestation
ph\'{e}nom\'{e}nale, mais Ren\'{e} Thom a ajout\'{e} \`{a} ce fondement de
l'exp\'{e}rience sensible l'intuition des discontinuit\'{e}s qualitatives
comme brisures de sym\'{e}trie de l'homog\'{e}n\'{e}it\'{e} spatiale et cas
g\'{e}n\'{e}ralis\'{e}s de transitions de phase (Thom [1972], [1980]).%
\footnote{%
Pour une introduction aux multiples facettes de la philosophie de Ren\'{e}
Thom, on pourra se r\'{e}f\'{e}rer aux actes du Colloque d'hommage tenu
\`{a} Cerisy en 1982 ainsi qu'\`{a} leur pr\'{e}sentation synth\'{e}tique
dans \textit{Physique du Sens} (Petitot [1992]).}

Dans mon hommage \`{a} Ren\'{e} Thom de 1994 dans l'ouvrage \textit{Passions
des formes} \'{e}dit\'{e} par Mich\`{e}le Porte, j'insistais beaucoup sur
l'actualit\'{e} scientifique des th\`{e}ses thomiennes dans les sciences
cognitives, en particulier en ce qui concerne les structures de la vision.
J'aimerais aujourd'hui dans cette courte note commenter un exemple allant
non pas de la th\'{e}orie des singularit\'{e}s vers la th\'{e}orie de la
vision mais de la th\'{e}orie de la vision vers la th\'{e}orie des
singularit\'{e}s.

En effet, comme pour toutes les id\'{e}es s\'{e}minales de Ren\'{e} Thom,
l'approfon\-dissement du c\^{o}t\'{e} perceptif de la th\'{e}orie
morphologique a conduit \`{a} des probl\`{e}mes math\'{e}matiques originaux.
Pour comprendre le r\^{o}le des singularit\'{e}s dans la perception des
formes, il a fallu faire la synth\`{e}se entre la classification de
Thom-Mather et le point de vue naturel de l'analyse g\'{e}om\'{e}trique des
images en vision computationnelle, qui est celui des mod\`{e}les dits 
\textit{multi-\'{e}chelle}. La synth\`{e}se a \'{e}t\'{e} effectu\'{e}e par
James Damon dans une s\'{e}rie d'articles et en particulier dans le papier
``Local Morse theory for solutions to the heat equation and Gaussian
blurring'' (Damon, 1995) sur lequel nous nous appuyerons ici.

La classification des mod\`{e}les multi\'{e}chelle de singularit\'{e}s fait
appara\^{i}tre certains ph\'{e}nom\`{e}nes nouveaux dont le plus surprenant
est celui de la possibilit\'{e} de cr\'{e}ation de points critiques lorsque
l'\'{e}chelle augmente, \'{e}v\'{e}nement qui semble violer le principe du
maximum auquel satisfont ces mod\`{e}les.\ D'o\`{u} l'int\'{e}r\^{e}t de
pr\'{e}ciser ce ph\'{e}nom\`{e}ne curieux.\ Nous le ferons ici au moyen de
graphiques calcul\'{e}s avec \textit{Mathematica}$^{\text{\textsc{tm}}}$.

Ces petits calculs ont \'{e}t\'{e} effectu\'{e}s il y a quelques ann\'{e}es
lors d'un colloque \`{a}, la Fondation des Treilles organis\'{e} avec
Bernard Teissier, Jean-Michel Morel et David Mumford dans le cadre du
Trimestre sp\'{e}cial du Centre Emile Borel ``Questions math\'{e}matiques en
traitement du signal et de l'image'' (septembre-d\'{e}cembre 1998). Ils
devaient servir de base \`{a} des discussions avec Ren\'{e} Thom mais la vie
en a d\'{e}cid\'{e} autrement.

\section{Le concept de ``Scale space analysis''}

L'id\'{e}e d'analyse multi\'{e}chelle, de ``scale space analysis'' ou
d'algorithmes ``pyramidaux'', domine les th\'{e}ories de l'analyse
g\'{e}om\'{e}trique des images. Elle remonte \`{a} Witkin (1983) et \`{a}
Koenderink (1984, 1986) et part de la remarque suivante. Pour \^{e}tre
morphologiquement correcte, une analyse des images doit s'effectuer en
termes de g\'{e}om\'{e}trie diff\'{e}rentielle. Le probl\`{e}me est que les
outils de cette derni\`{e}re ne sont pas directement applicables au signal
en tant que tel, qui est trop bruit\'{e} pour \^{e}tre diff\'{e}rentiable.
Pour que les images puissent acqu\'{e}rir le statut d'observables
g\'{e}om\'{e}triquement analysables par d\'{e}tection d'invariants, il faut
par cons\'{e}quent d\'{e}finir au pr\'{e}alable une \'{e}chelle, c.\`{a}.d.
fixer un niveau de r\'{e}gularisation du signal bruit\'{e}. Jan Koenderink
l'a souvent soulign\'{e}~: le scaling est essentiel. C'est ce que fait
d'ailleurs la vision naturelle puisque les neurones visuels primaires
poss\`{e}dent un champ r\'{e}cepteur qui est un petit domaine de
photor\'{e}cepteurs r\'{e}tiniens et op\`{e}rent sur le signal optique comme
des filtres par convolution avec leur profil r\'{e}cepteur (leur fonction de
transfert) qui est une fonction d\'{e}finie sur leur champ r\'{e}cepteur.

Mais comment une analyse multi\'{e}chelle peut-elle d\'{e}boucher sur une
v\'{e}ritable analyse morphog\'{e}n\'{e}tique d'une image 2D d\'{e}finie
comme une fonction (ou, mieux, une distribution) $I(x,y)$ sur la fen\^{e}tre
r\'{e}tinienne $R$ de coordonn\'{e}es $x$ et $y$? L'id\'{e}e directrice est
de plonger l'image dans une famille $I_{s}(x,y)$ param\'{e}tr\'{e}e par une
\'{e}chelle $s$ de fa\c{c}on \`{a} ce que~:

\begin{description}
\item[(i)]  $I_{0}=I$,

\item[(ii)]  $I_{1}$ soit une image indiff\'{e}renci\'{e}e, et

\item[(iii)]  Lorsque l'\'{e}chelle $s$ cro\^{\i }t, l'\'{e}volution de $%
I_{0}$ \`{a} $I_{1}$ ``simplifie'' strictement l'image.\ Cette contrainte
dite de ``causalit\'{e}'' interdit l'apparition \textit{ex nihilo} de
nouveaux d\'{e}tails lorsque l'\'{e}chelle cro\^{\i }t.
\end{description}

L'\'{e}volution avec l'\'{e}chelle $s$ des lignes de niveau de $I_{s}$,
c'est-\`{a}-dire la suite d'\'{e}v\'{e}nements de bifurcation qu'elles
subissent en se simplifiant progressivement, fournit une m\'{e}thode
puissante pour analyser la structure morphologique de l'image et sa
d\'{e}composition en \'{e}l\'{e}ments constituants. On montre que, sous des
contraintes g\'{e}n\'{e}rales de lin\'{e}arit\'{e}, d'invariance par
translation, d'isotropie et d'invariance d'\'{e}chelle, la fa\c{c}on la plus
simple d'obtenir un tel r\'{e}sultat est de prendre pour $I_{s}$ une
solution de \textit{l'\'{e}quation de diffusion} typique qu'est
l'\'{e}quation de la chaleur $\partial I_{s}/\partial s=\Delta I_{s}$. Dans
la mesure o\`{u} le noyau de la chaleur est gaussien, on est ainsi conduit
\`{a} l'id\'{e}e d'un lissage gaussien multi\'{e}chelle de l'image (Gaussian
blurring).\footnote{%
Il faut insister sur le fait que $s$ est ici un param\`{e}tre d'\'{e}chelle
et non pas un param\`{e}tre temporel comme c'est le cas d'habitude.\ L'\'{e}%
volution s'effectue dans un espace-\'{e}chelle et non pas dans un
espace-temps.}

Mais ce qui domine ph\'{e}nom\'{e}nologiquement une image sont les
singularit\'{e}s.\ Si l'on veut adapter \`{a} l'analyse d'images la
th\'{e}orie de Thom-Mather des singularit\'{e}s g\'{e}n\'{e}riques et des
d\'{e}ploiements universels, il faut par cons\'{e}quent la rendre
multi\'{e}chelle. 

\section{Th\'{e}orie multi\'{e}chelle des singularit\'{e}s}

James Damon a montr\'{e} comment on pouvait transformer la th\'{e}orie de
Morse-Whitney-Thom-Mather-Arnold dans le cas de l'\'{e}quation de la
chaleur. La difficult\'{e} principale est que les formes normales de Morse
ne satisfont pas \`{a} cette EDP. Les m\'{e}thodes doivent donc \^{e}tre
transpos\'{e}es des espaces de germes d'applications $C^{\infty }$ aux
espaces de germes de solutions de l'EDP, or ces espaces n'ont pas les
``bonnes'' propri\'{e}t\'{e}s alg\'{e}briques qui font marcher la
th\'{e}orie et permettent d'appliquer des th\'{e}or\`{e}mes de
transversalit\'{e} de Thom pour obtenir des r\'{e}sultats de
g\'{e}n\'{e}ricit\'{e}.\footnote{%
Pour une premi\`{e}re introduction \`{a} la th\'{e}orie des singularit\'{e}s
de Thom-Mather on pourra consulter en particulier le s\'{e}minaire Bourbaki
Chenciner [1973], l'article de l'\textit{Encyclop\ae dia Universalis }%
Chenciner [1980] et le Graduate Text Golubitsky\&Guillemin [1973]. Les six
articles de base de John Mather sont inclus dans la bibliographie.}

James Damon a donc d\^{u} d'abord red\'{e}finir le concept d'\'{e}quivalence
pour des germes $C^{\infty }$ d'applications $f(x,s)$

\[
\begin{array}{cccc}
f: & \left( \Bbb{R}^{n+1},0\right) & \longrightarrow & \left( \Bbb{R}%
,0\right) \\ 
& (x,s) & \longmapsto & z=f(x,s)
\end{array}
\]

\noindent o\`{u} $x=\left( x_{1},\ldots ,x_{n}\right) .$

\noindent $\bullet $ $f$\textit{\ }et $g$\textit{\ }sont $H$%
-\'{e}quivalentes s'il existe un germe de diff\'{e}omorphisme $\varphi
:\left( \Bbb{R}^{n+1},0\right) \rightarrow \left( \Bbb{R}^{n+1},0\right) $de
la forme $\varphi (x,s)=\left( \varphi _{1}(x,s),\varphi _{2}(s)\right) $
avec $\varphi _{2}^{\prime }(0)>0$ et un germe $c:\left( \Bbb{R},0\right)
\rightarrow \left( \Bbb{R},0\right) $ tels que 
\[
g(x,s)=f\circ \varphi (x,s)+c(s). 
\]

\noindent Autrement dit, la diff\'{e}rence entre $g$ et $f$ peut \^{e}tre
r\'{e}sorb\'{e}e

\begin{description}
\item[(i)]  par des translations au but $c(s)$ d\'{e}pendant de $s$,

\item[(ii)]  par un changement de coordonn\'{e}e (pr\'{e}servant
l'orientation) $\varphi _{2}$ sur l'axe $s$ de la source,

\item[(iii)]  par des changements de coordonn\'{e}es $\varphi _{1}(\bullet %
,s),$ d\'{e}pendant de $s$, du sous-espace $\Bbb{R}^{n}$ de la source.
\end{description}

\noindent $\bullet $ De m\^{e}me $f$ et $g$ sont $IS$-\'{e}quivalentes ($IS$
signifie ``Intensity Sensitive'') si $c$ est constant~\footnote{$c=0$ si $%
g(0,0)=0$, mais on peut avoir \`{a} consid\'{e}rer des cas $c\neq 0.$} et
s'il existe de plus un germe de diff\'{e}omorphisme $\psi :\left( \Bbb{R}%
^{2},0\right) \rightarrow \left( \Bbb{R}^{2},0\right) $ de la forme $\psi
(z,s)=\left( \theta (z,s),s\right) $ avec $\frac{\partial \theta }{\partial z%
}(0,0)>0$ et $\theta (0,s)=0$ pour tout $s$ tels que

\[
g(x,s)=\theta \circ f\circ \varphi (x,s)+c. 
\]

\noindent Cela signifie que l'on peut remplacer $c(s)$ par des changements
de coordonn\'{e}es appropri\'{e}s au but (pr\'{e}servant l'orientation) $%
\theta (\bullet ,s)$ d\'{e}pendant de $s$.

\`{A} ces notions d'\'{e}quivalence sont naturellement associ\'{e}es des
notions de stabilit\'{e} par d\'{e}formation.\ On consid\`{e}re des familles 
$f_{w}(x,s)$ de fonctions $f(x,s)$ param\'{e}tr\'{e}es par des
param\`{e}tres $w\in \Bbb{R}^{q}$

\[
\begin{array}{cccc}
f_{w}: & \left( \Bbb{R}^{n+1+q},0\right) & \longrightarrow & \left( \Bbb{R}%
,0\right) \\ 
& (x,s,w) & \longmapsto & z=f(x,s,w)=f_{w}(x,s)
\end{array}
\]

\noindent et d\'{e}formant $f(x,s,0)=f_{0}(x,s)$.

\textit{D\'{e}finition}. $f_{0}(x,s)$ est $H$-stable (resp. $IS$-stable) si $%
f_{0}(x,s)$ est son propre d\'{e}ploiement universel, autrement dit si toute
d\'{e}formation assez petite est triviale pour la $H$-\'{e}quivalence (resp.
la $IS$-\'{e}quivalence).

En adaptant les d\'{e}monstrations de Ren\'{e} Thom et de John Mather, James
Damon a trouv\'{e} les formes normales des singularit\'{e}s stables et des
d\'{e}ploiements universels de bas degr\'{e}. Par exemple les
\'{e}quivalents $H$-stables des points critiques quadratiques (points
critiques non d\'{e}g\'{e}n\'{e}r\'{e}s stables ) pour la $H$%
-\'{e}quivalence ont pour forme normale~:

\[
\pm \sum\limits_{i=1}^{i=n}x_{i}^{2}\pm (2n)s 
\]

\noindent et

\[
\sum\limits_{i=1}^{i=n}a_{i}x_{i}^{2}\text{ avec }a_{i}\neq 0\text{ pour
tout }i\text{ et }\sum\limits_{i=1}^{i=n}a_{i}=0 
\]

\noindent dont on v\'{e}rifie trivialement qu'ils satisfont l'\'{e}quation
de la chaleur. Mais on obtient aussi des singularit\'{e}s $H$-stables 
\textit{cubiques} de forme

\[
x_{1}^{3}+6sx_{1}+Q(x_{2},\ldots ,x_{n},s) 
\]

\noindent et

\[
x_{1}^{3}-6sx_{1}+-6x_{1}x_{2}^{2}+Q(x_{2},\ldots ,x_{n},s) 
\]

\noindent o\`{u} $Q(x_{2},\ldots ,x_{n},s)$ est une singularit\'{e}
quadratique en $(x_{2},\ldots ,x_{n}).\;$Cela est d\^{u} au fait que la
variable d'\'{e}chelle $s$ est de poids $2$ relativement aux variables $%
x_{i} $ qui sont de poids $1$ et qu'un terme apparemment quadratique comme $%
sx_{i}$ est donc en fait cubique.

Pour la $IS$-\'{e}quivalence, on obtient la m\^{e}me premi\`{e}re forme de
singularit\'{e}s quadratiques stables

\[
\pm \sum\limits_{i=1}^{i=n}x_{i}^{2}\pm (2n)s 
\]

\noindent Mais la seconde forme n'est plus la m\^{e}me et devient

\[
\sum\limits_{i=1}^{i=n}a_{i}x_{i}^{2}+2\left(
\sum\limits_{i=1}^{i=n}a_{i}\right) s\text{ avec }a_{i}\neq 0\text{ pour
tout }i\text{ et }\sum\limits_{i=1}^{i=n}a_{i}\neq 0 
\]

\noindent Pour $\sum\limits_{i=1}^{i=n}a_{i}=0$ on obtient une nouvelle
forme plus compliqu\'{e}e qui, dans le cas $n=2$, s'\'{e}crit~:

\[
\sum\limits_{i=1}^{i=2}a_{i}x_{i}^{2}\pm \left( s^{2}+\frac{1}{2}s\left(
\sum\limits_{i=1}^{i=2}x_{i}^{2}\right) +\frac{1}{16}\left(
\sum\limits_{i=1}^{i=2}x_{i}^{2}\right) ^{2}\right) 
\]

\noindent  $\text{ avec }a_{i}\neq 0%
\text{ pour tout }i\text{ et }\sum\limits_{i=1}^{i=n}a_{i}=0.$  On obtient \'{e}galement deux formes cubiques $IS$-stables~:

\[
x_{1}^{3}+6sx_{1}+\sum\limits_{i=2}^{i=n}a_{i}x_{i}^{2}+2\left(
\sum\limits_{i=2}^{i=n}a_{i}\right) s\text{ avec }a_{i}\neq 0\text{ pour
tout }i\text{ et }\sum\limits_{i=2}^{i=n}a_{i}\neq 0 
\]

\noindent et

\[
x_{1}^{3}-6sx_{1}+-6x_{1}x_{2}^{2}+\sum\limits_{i=2}^{i=n}a_{i}x_{i}^{2}+2%
\left( \sum\limits_{i=2}^{i=n}a_{i}\right) s\text{ avec }a_{i}\neq 0\text{
pour tout }i\text{ et }\sum\limits_{i=2}^{i=n}a_{i}\neq 0 
\]

C'est cette derni\`{e}re forme qui est surprenante car elle conduit \`{a}
des \textit{cr\'{e}ations} de points critiques, ce qui pourrait sembler
interdit par le principe de causalit\'{e} expos\'{e} plus haut.

\section{Analyse graphique de la cr\'{e}ation de points critiques par
diffusion}

Nous allons analyser graphiquement la cr\'{e}ation de points critiques dans
le cas $n=2$. Nous \'{e}tudions donc la solution de l'\'{e}quation de la
chaleur~:

\[
f(x,y,s)=x^{3}-6xy^{2}+y^{2}-6sx+2s 
\]

\noindent consid\'{e}r\'{e}e comme une famille \`{a} un param\`{e}tre
(l'\'{e}chelle $s$) de fonctions de $(x,y)$.

\subsection{Forme des points critiques}

Les points critiques sont donn\'{e}s par 
\[
\left\{ 
\begin{array}{l}
\frac{\partial f}{\partial x}=3x^{2}-6y^{2}-6s=0 \\ 
\frac{\partial f}{\partial y}=-12xy+2y=0
\end{array}
\right. 
\]

\noindent dont les $4$ solutions sont~: 
\[
\left\{ 
\begin{array}{l}
x_{1\pm }=\frac{1}{6}\text{, }y_{1\pm }=\pm \frac{(1-72s)^{1/2}}{6\sqrt{2}}
\\ 
x_{2\pm }=\pm \sqrt{2s}\text{, }y_{2\pm }=0
\end{array}
\right. 
\]

\begin{figure}[tbp]
\begin{center}
\includegraphics[
width=1.8265in,height=2.0141in]%
{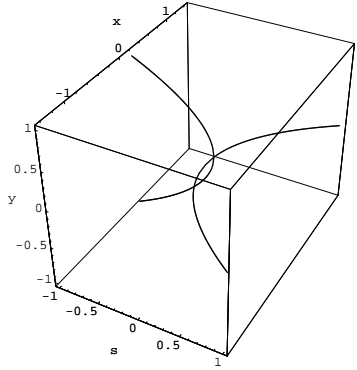}
\end{center}
\caption{Les points critiques de $f(x,y,s)$ se regroupent en $2$ branches.}
\label{BranchesPC}
\end{figure}

\noindent Ils constituent deux branches paraboliques param\'{e}tris\'{e}es
par $s$, $pc_{1}(s)$ et $pc_{2}(s)$.\ Pour $s<0$ (ce qui est interdit pour
une \'{e}chelle), $pc_{2}(s)$ est imaginaire.\ Pour $s>1/72$ c'est la
branche $pc_{1}$ qui devient imaginaire, $pc_{2}$ donnant deux points
critiques oppos\'{e}s.\ Pour $0\leq s\leq 1/72$, il y a $2$ paires de points
critiques (index\'{e}es par $1\pm $ et $2\pm $)~: la branche $pc_{2}$
appara\^{i}t (cr\'{e}ation de deux points critiques oppos\'{e}s) et reste en
comp\'{e}tition avec $pc_{1}$ tant que $s<1/72$. 
Cf. Figure \ref{BranchesPC}.

Le Hessien $H$ de $f$ \'{e}tant donn\'{e} par 
\[
H=\left( \frac{\partial ^{2}f}{\partial x_{i}\partial x_{i}}\right) =\left( 
\begin{array}{ll}
6x & -12y \\ 
-12y & 2-12x
\end{array}
\right) , 
\]

\noindent sa valeur aux points critiques est 
\[
H_{1\pm }=\left( 
\begin{array}{ll}
1 & \mp \sqrt{2(1-72s)} \\ 
\mp \sqrt{2(1-72s)} & 0
\end{array}
\right) \text{ et }H_{2\pm }=\left( 
\begin{array}{ll}
\pm 6\sqrt{2s} & 0 \\ 
0 & 2\left( 1\mp 6\sqrt{2s}\right)
\end{array}
\right) 
\]

\noindent et ses valeurs propres (dans l'ordre de $H_{1+}$, $H_{1-}$, $%
H_{2+} $, $H_{2-}$) sont

\begin{eqnarray*}
&&\left\{ \frac{1-3\sqrt{1-64s}}{2},\frac{1+3\sqrt{1-64s}}{2}\right\}
,\left\{ \frac{1+3\sqrt{1-64s}}{2},\frac{1-3\sqrt{1-64s}}{2}\right\} , \\
&&\left\{ 2-12\sqrt{2s},6\sqrt{2s}\right\} ,\left\{ 2+12\sqrt{2s},-6\sqrt{2s}%
\right\} .
\end{eqnarray*}

\noindent Quant aux vecteurs propres, on trouve, toujours pour le m\^{e}me
ordre,

\begin{eqnarray*}
&&\left\{ \left\{ \frac{-1+3\sqrt{1-64s}}{2\sqrt{2(1-72s)}},1\right\}
,\left\{ \frac{-1-3\sqrt{1-64s}}{2\sqrt{2(1-72s)}},1\right\} \right\} , \\
&&\left\{ \left\{ \frac{1-3\sqrt{1-64s}}{2\sqrt{2(1-72s)}},1\right\}
,\left\{ \frac{1+3\sqrt{1-64s}}{2\sqrt{2(1-72s)}},1\right\} \right\} , \\
&&\left\{ \left\{ 0,1\right\} ,\left\{ 1,0\right\} \right\} ,\left\{ \left\{
0,1\right\} ,\left\{ 1,0\right\} \right\} .
\end{eqnarray*}

Comme les deux points critiques de la branche $pc_{1}$ deviennent
imaginaires conjugu\'{e}s pour $s>1/72,$ ce sera a fortiori le cas pour $%
s>1/64$. Les changements de signe des valeurs propres des points critiques
r\'{e}els s'effectuent pour $1-3\sqrt{1-64s}=0$ et $1-6\sqrt{2s}=0$, ce qui
dans les deux cas donne $s=1/72$. D'o\`{u} les tableaux suivants o\`{u} $%
\lambda $ et $\mu $ sont les deux valeurs propres.

\begin{description}
\item[(i)]  Pour $0<s<1/72$~: 
\[
\begin{array}{cccc}
\begin{array}{c}
pc_{1+}= \\ 
\left( \frac{1}{6},\frac{(1-72s)^{1/2}}{6\sqrt{2}}\right) 
\end{array}
& 
\begin{array}{c}
pc_{1-}= \\ 
\left( \frac{1}{6},-\frac{(1-72s)^{1/2}}{6\sqrt{2}}\right) 
\end{array}
& 
\begin{array}{c}
pc_{2+}= \\ 
\left( \sqrt{2s},0\right) 
\end{array}
& 
\begin{array}{c}
pc_{2-}= \\ 
-\left( \sqrt{2s},0\right) 
\end{array}
\\ 
\lambda <0 & \lambda >0 & \lambda >0 & \lambda >0 \\ 
\mu >0 & \mu <0 & \mu >0 & \mu <0 \\ 
\text{col} & \text{col} & \text{sommet} & \text{col}
\end{array}
\]

\item[(ii)]  Pour $s>1/72$~:
\end{description}

\[
\begin{array}{cc}
\begin{array}{c}
pc_{2+}= \\ 
\left( \sqrt{2s},0\right) 
\end{array}
& 
\begin{array}{c}
pc_{2-}= \\ 
-\left( \sqrt{2s},0\right) 
\end{array}
\\ 
\lambda <0 & \lambda >0 \\ 
\mu >0 & \mu <0 \\ 
\text{col} & \text{col}
\end{array}
\]

\subsection{Repr\'{e}sentations des surfaces}

Aux points critiques les valeurs de $f(x,y,s)$ sont $z_{1\pm }=s+1/216$, $%
z_{2\pm }=2s\left( 1\mp 2\sqrt{2s}\right) $. Pour bien visualiser la famille
de fonctions $f(x,y,s)$ nous allons repr\'{e}senter de plusieurs
mani\`{e}res les graphes de leurs sections \`{a} 
$s=\text{cste}$.

La planche $1$ represente $6$ sections pour $s$ variant de $-(1/3)1/72$
\`{a} ($4/3)1/72$ par pas de $(1/3)1/72$. Le point de vue par le dessus
permet de bien visualiser la position des points critiques.

\begin{tabular}{cc}

\leavevmode\lower0in\hbox{
\includegraphics[
width=1.6873in,height=2.0211in]%
{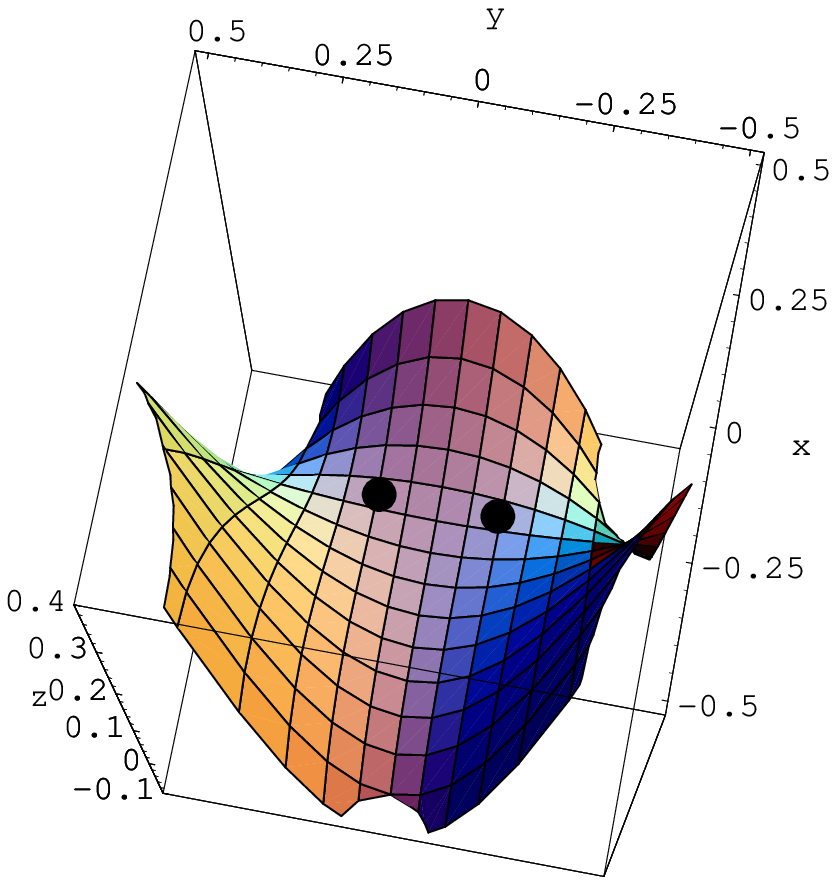}}

& 

\leavevmode\lower0in\hbox{
\includegraphics[
width=1.6873in,height=2.0211in]%
{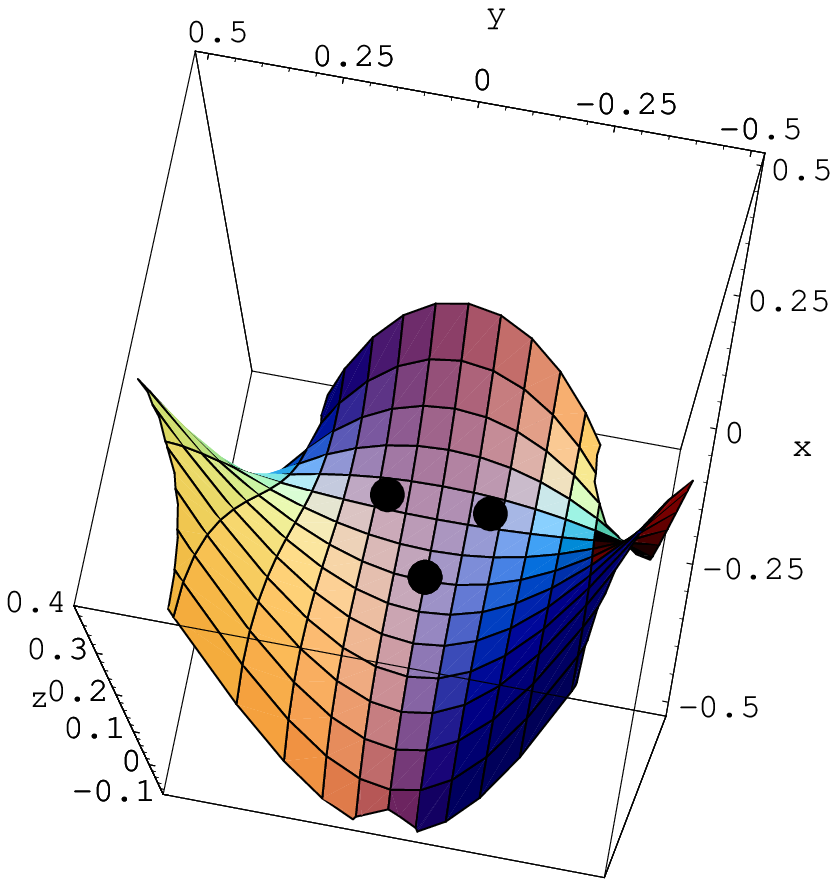}}

\\ 

\leavevmode\lower0in\hbox{
\includegraphics[
width=1.6873in,height=2.0211in]%
{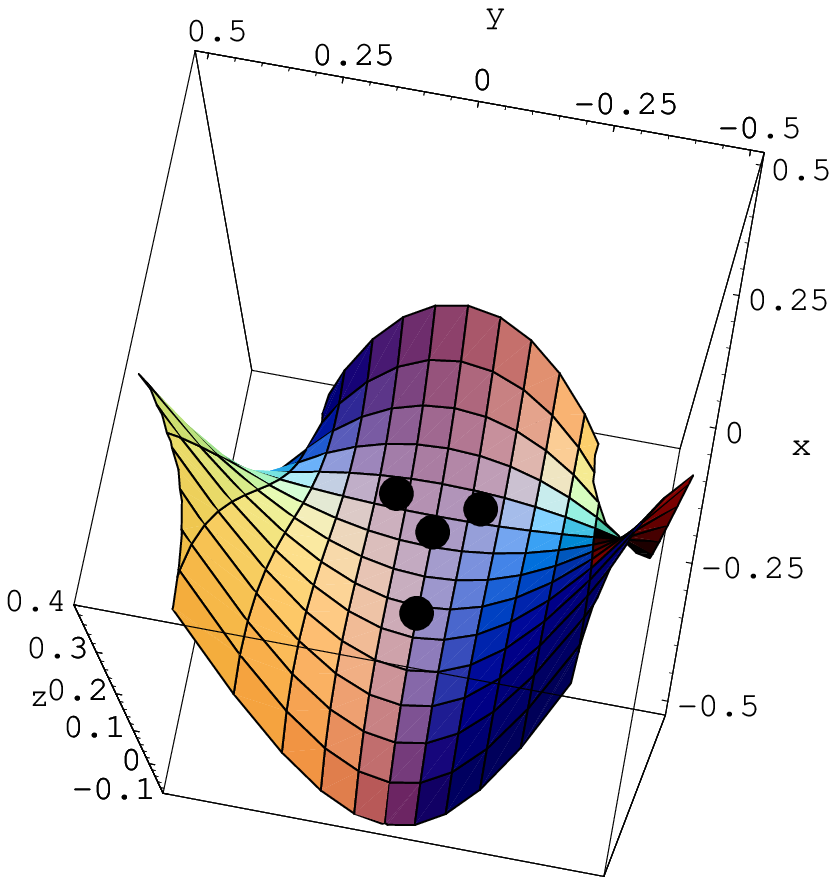}}

& 

\leavevmode\lower0in\hbox{
\includegraphics[
width=1.6873in,height=2.0211in]%
{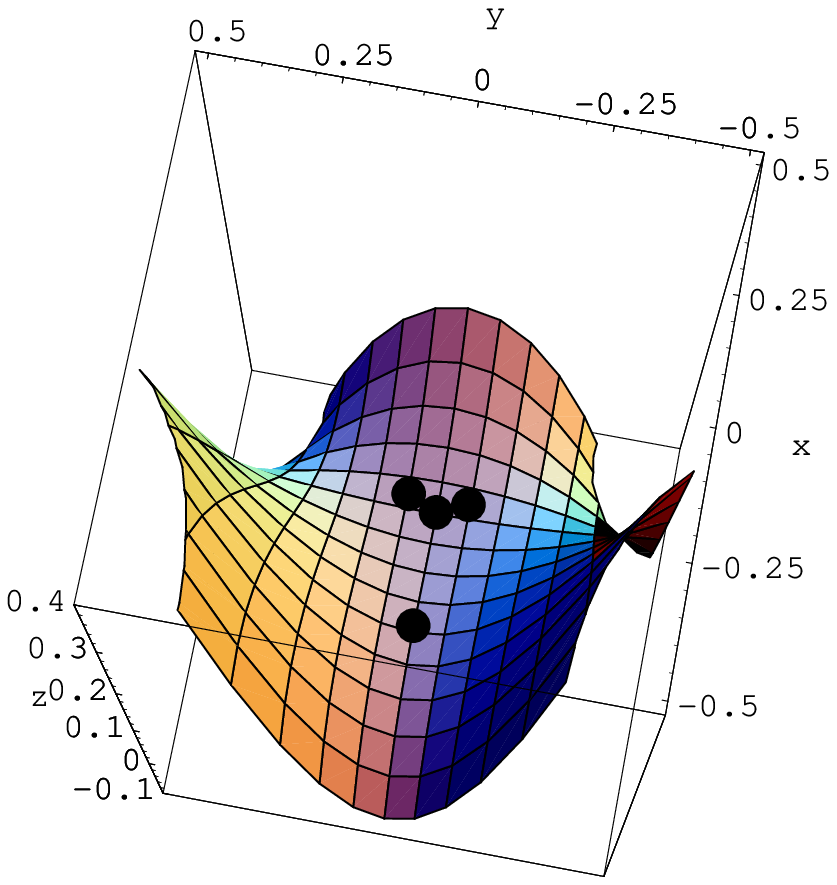}}

\\ 

\leavevmode\lower0in\hbox{
\includegraphics[
width=1.6873in,height=2.0211in]%
{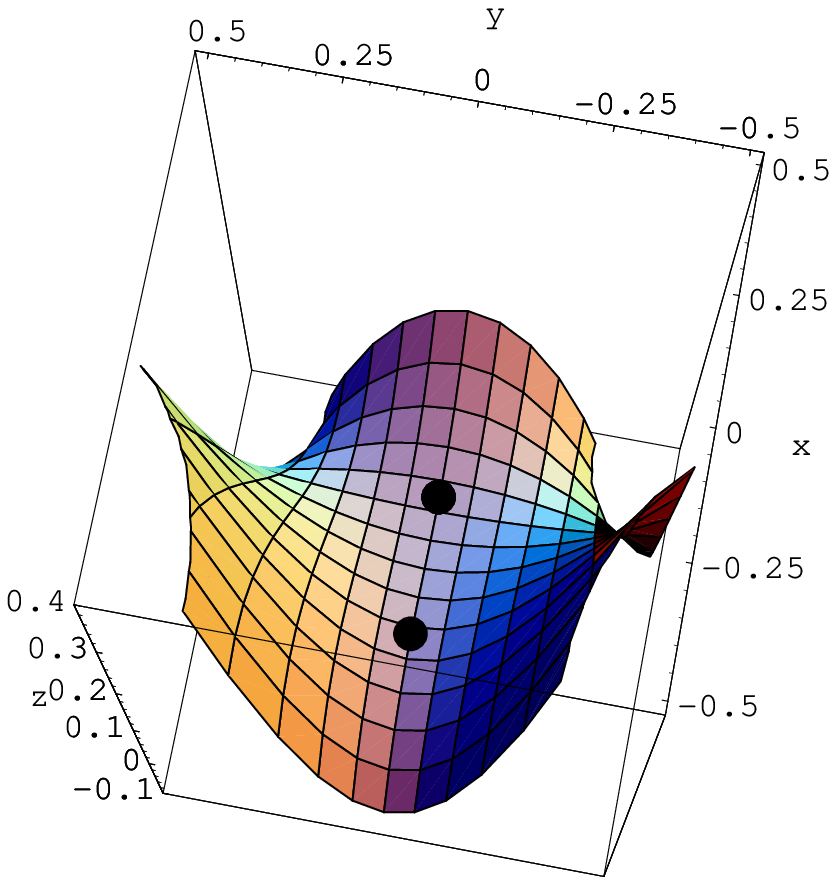}}

& 

\leavevmode\lower0in\hbox{
\includegraphics[
width=1.6873in,height=2.0211in]%
{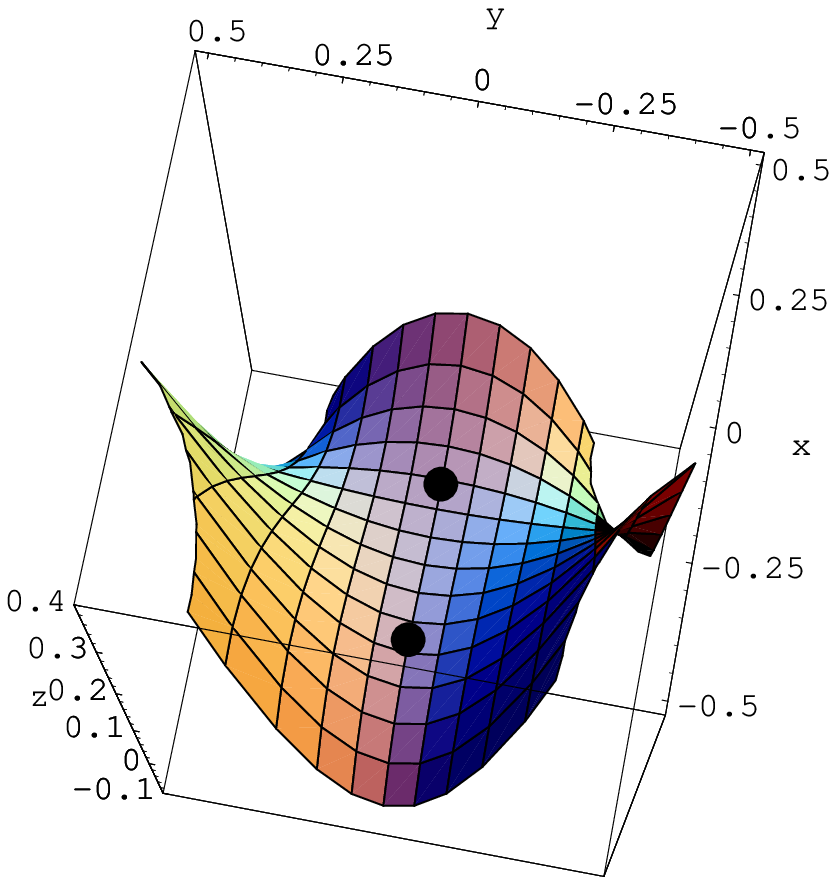}}

\\ 
\multicolumn{2}{c}{Planche 1}
\end{tabular}

La planche $2$ repr\'{e}sente les m\^{e}mes sections mais dans un point de
vue plus frontal permettant de mieux voir le relief de $f$.

\begin{tabular}{ccc}

\leavevmode\lower0in\hbox{
\includegraphics[
width=1.3474in,height=2.5348in]%
{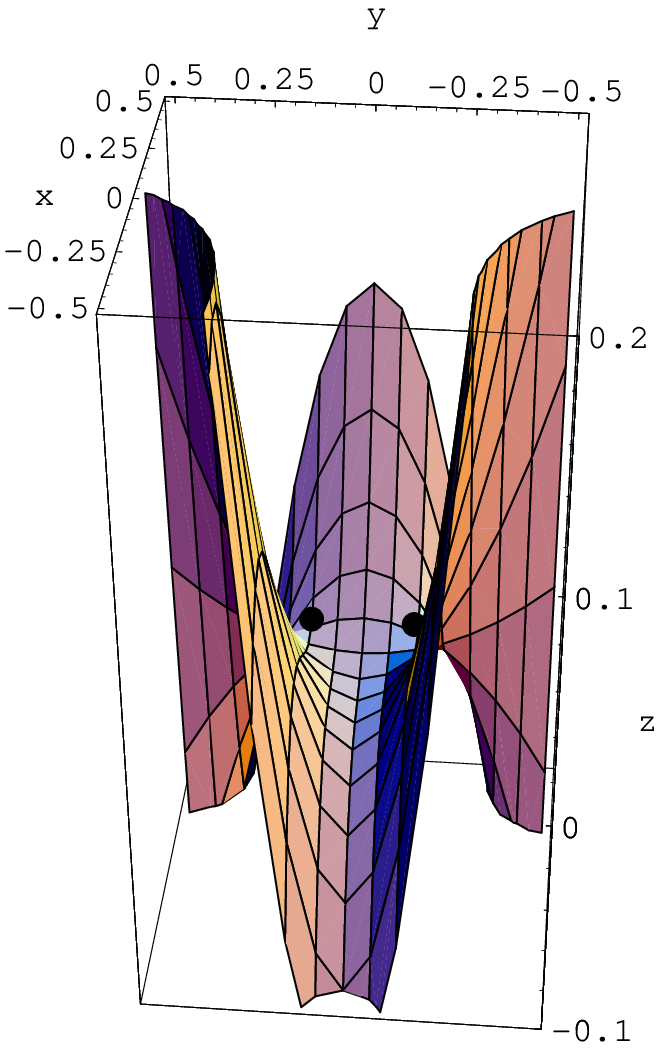}}

& 

\leavevmode\lower0in\hbox{
\includegraphics[
width=1.3474in,height=2.5348in]%
{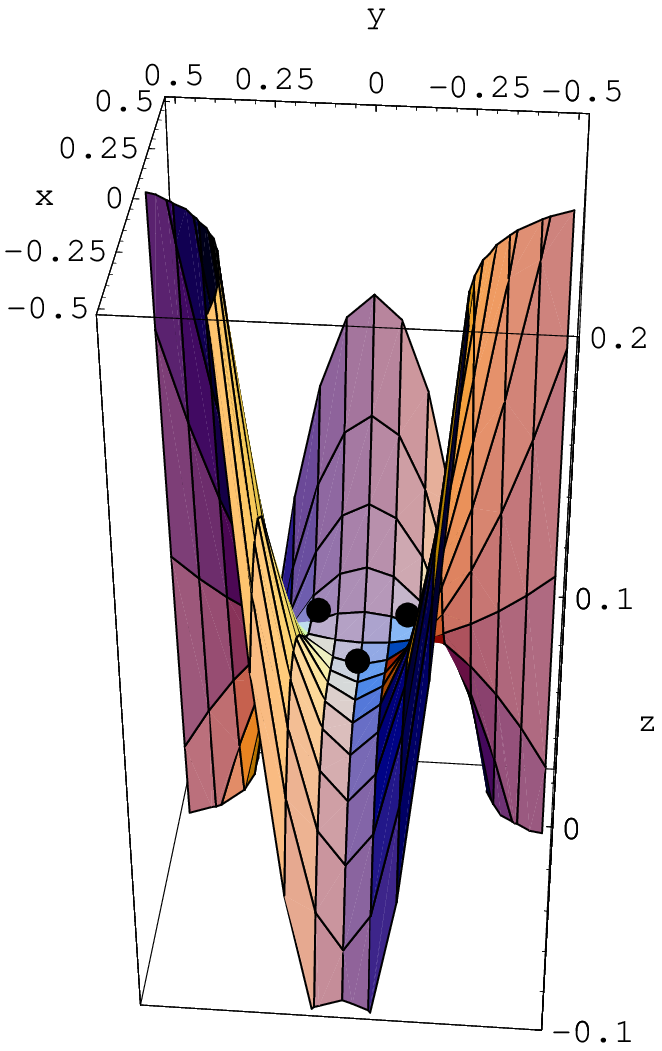}}

& 

\leavevmode\lower0in\hbox{
\includegraphics[
width=1.3474in,height=2.5348in]%
{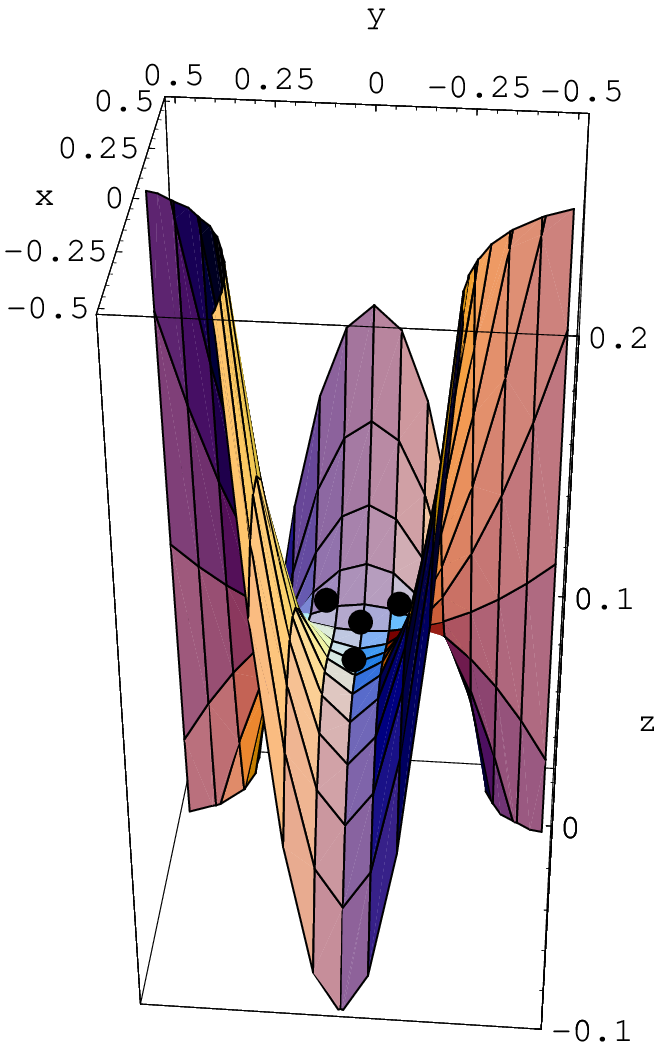}}

\\ 

\leavevmode\lower0in\hbox{
\includegraphics[
width=1.3474in,height=2.5348in]%
{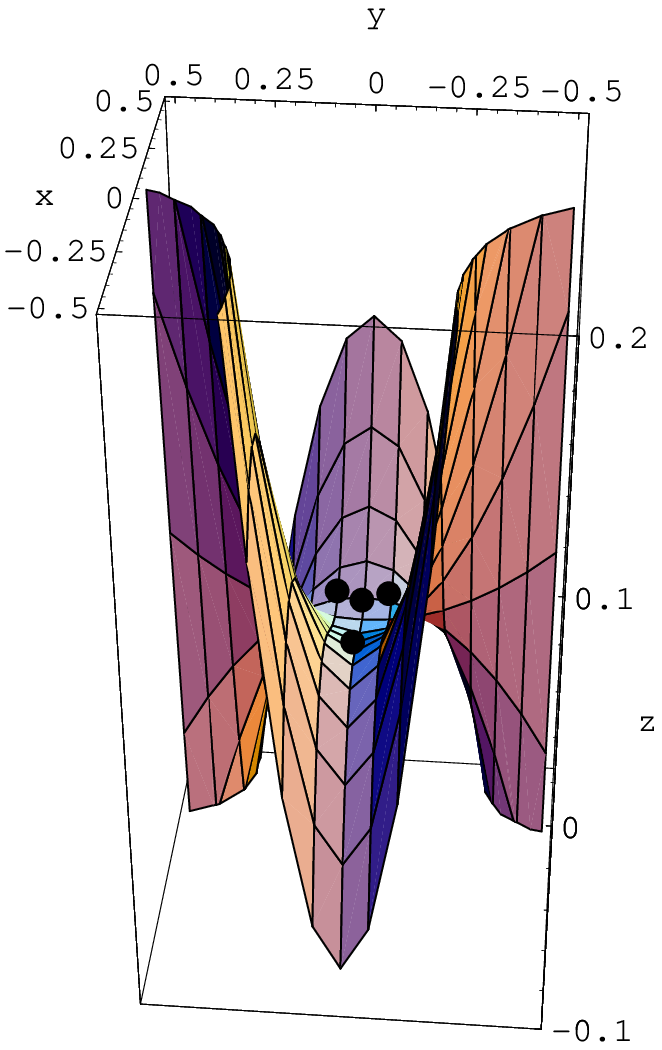}}

& 

\leavevmode\lower0in\hbox{
\includegraphics[
width=1.3474in,height=2.5348in]%
{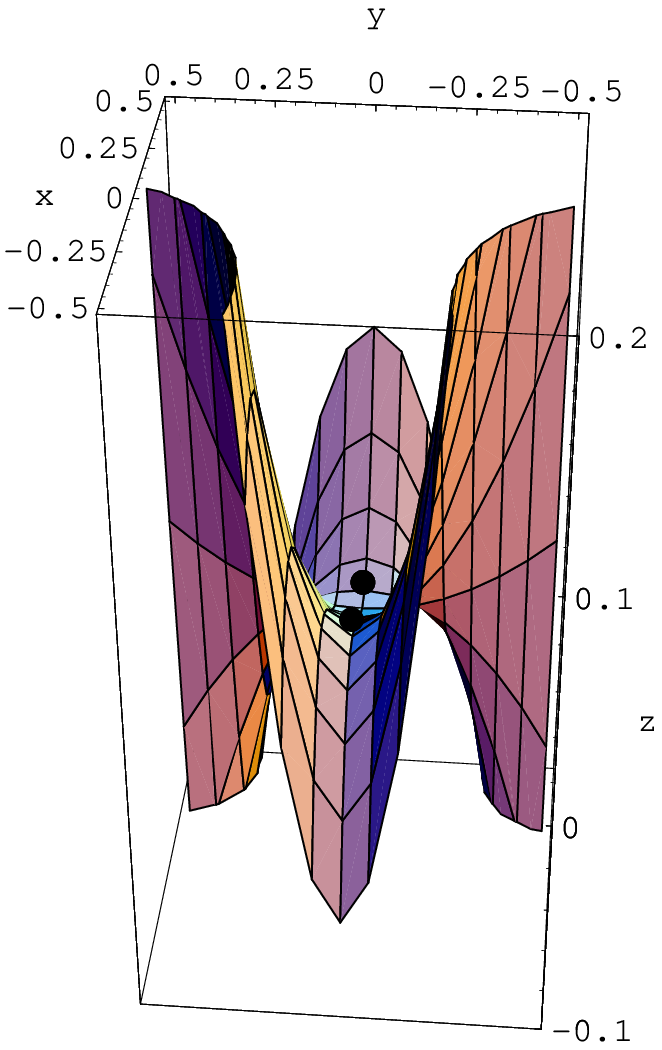}}

& 

\leavevmode\lower0in\hbox{
\includegraphics[
width=1.3474in,height=2.5348in]%
{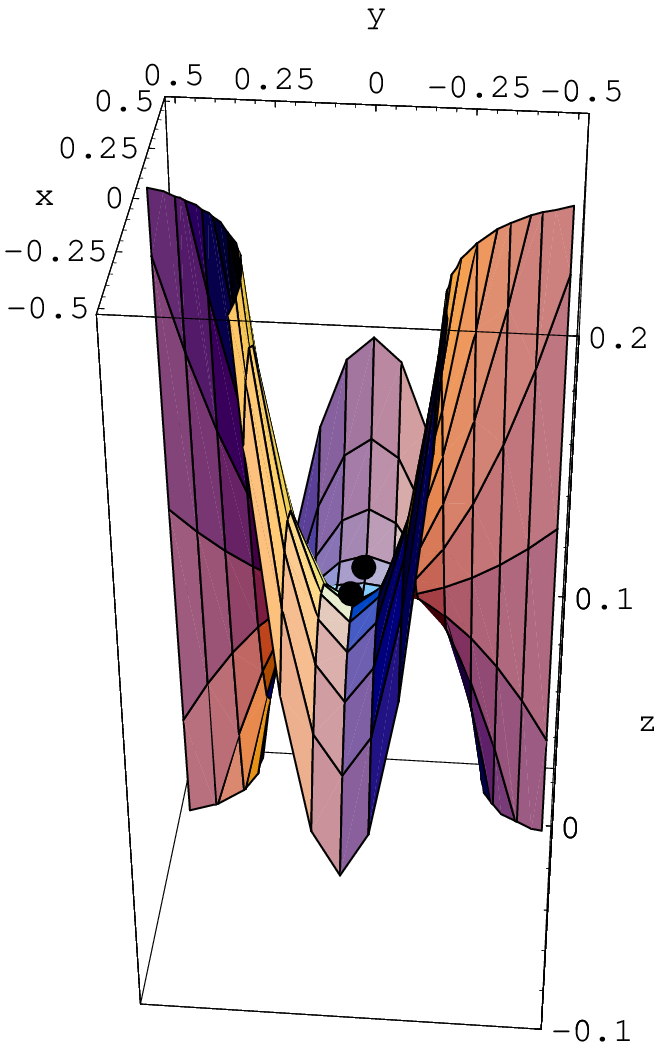}}

\\ 
\multicolumn{3}{c}{Planche 2}
\end{tabular}

Quant \`{a} la planche $3$ elle montre les lignes de niveaux de $f$ ainsi
que son champ de gradient $\nabla f$ par rapport aux variables spatiales $x$
et $y$.

\begin{tabular}{cc}

\leavevmode\lower0in\hbox{
\includegraphics[
width=2.0211in,height=2.0211in]%
{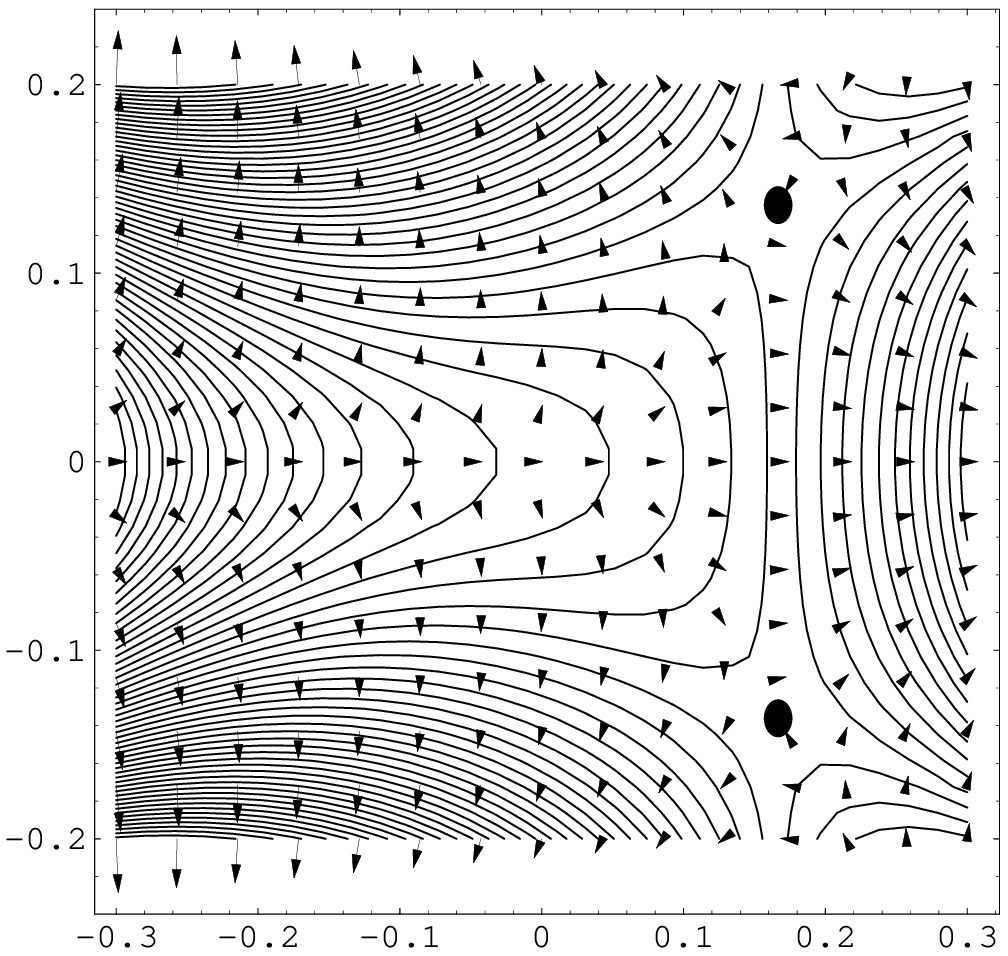}}

& 

\leavevmode\lower0in\hbox{
\includegraphics[
width=2.0211in,height=2.0211in]%
{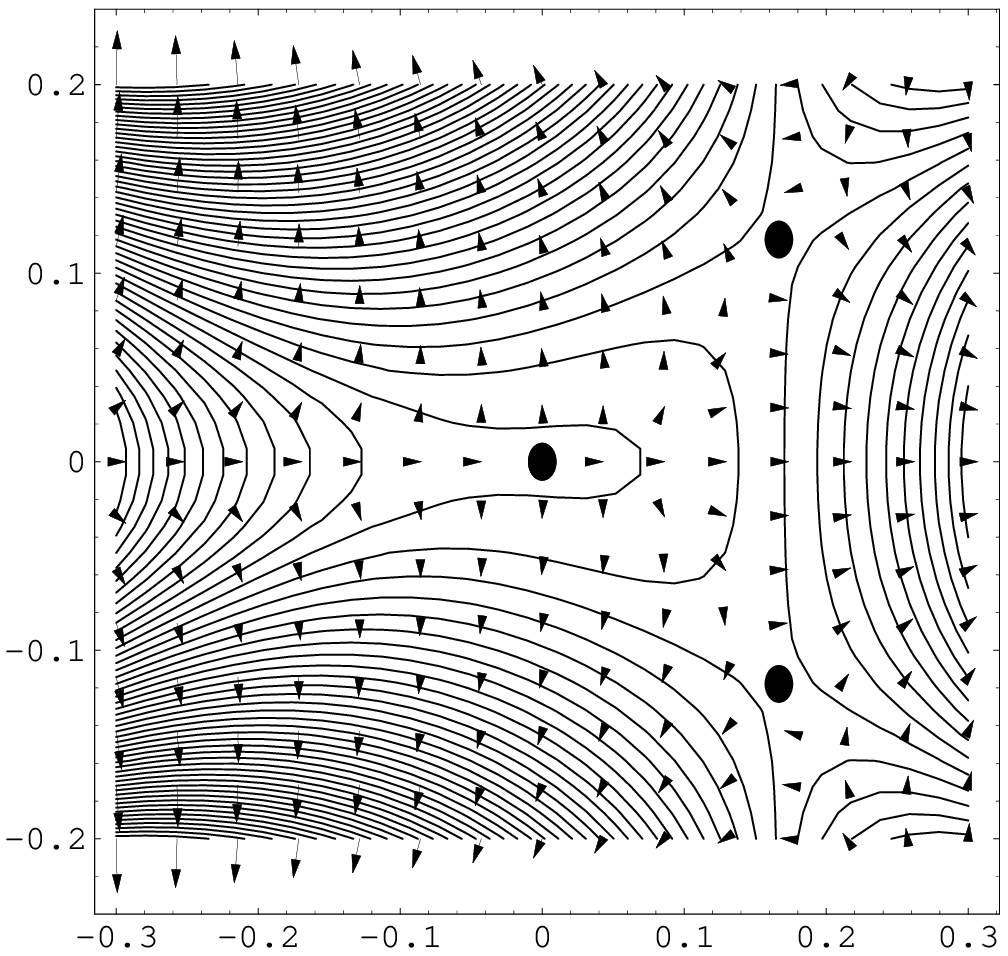}}

\\ 

\leavevmode\lower0in\hbox{
\includegraphics[
width=2.0211in,height=2.0211in]%
{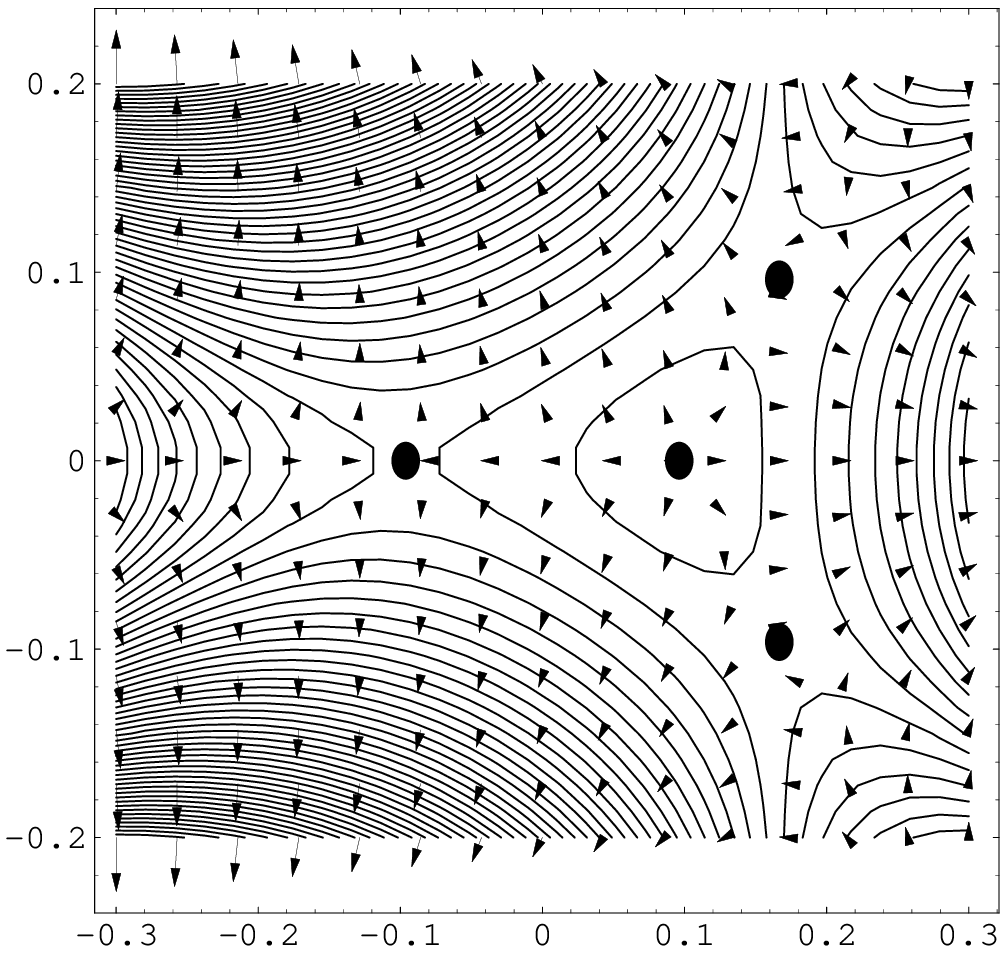}}

& 

\leavevmode\lower0in\hbox{
\includegraphics[
width=2.0211in,height=2.0211in]%
{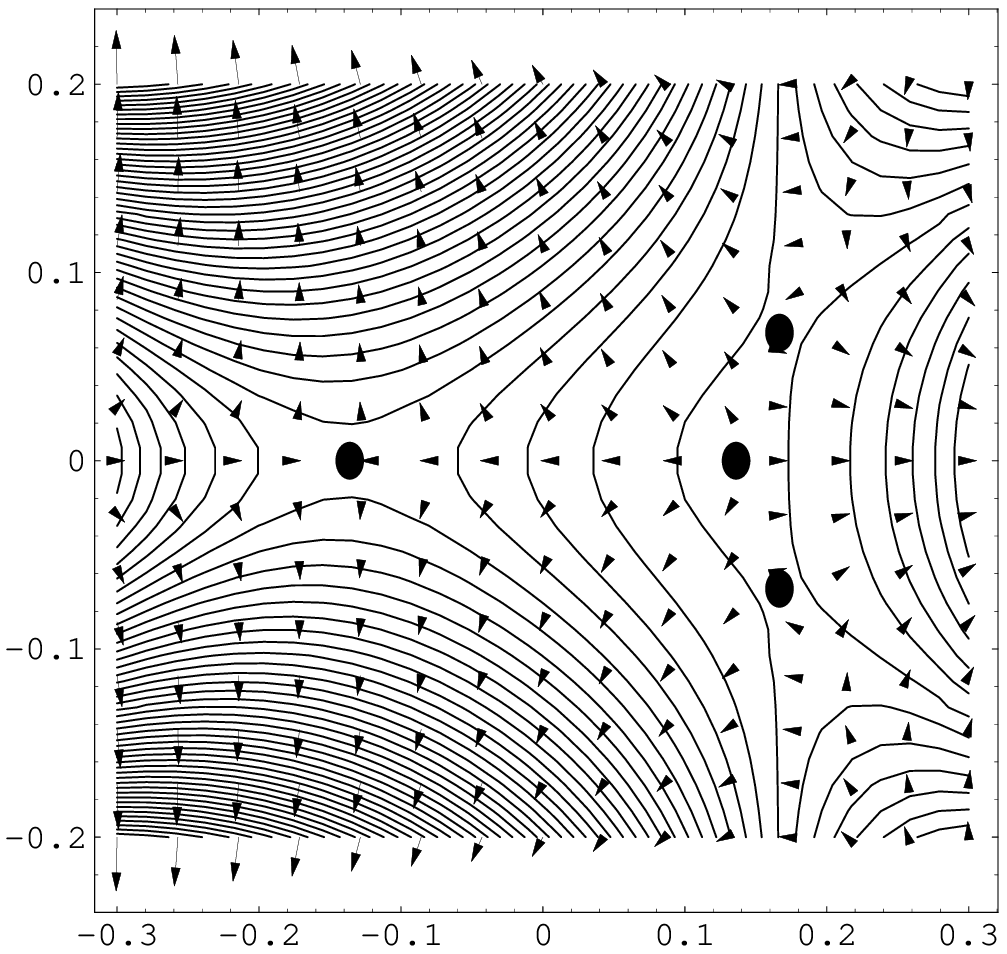}}

\\ 

\leavevmode\lower0in\hbox{
\includegraphics[
width=2.0211in,height=2.0211in]%
{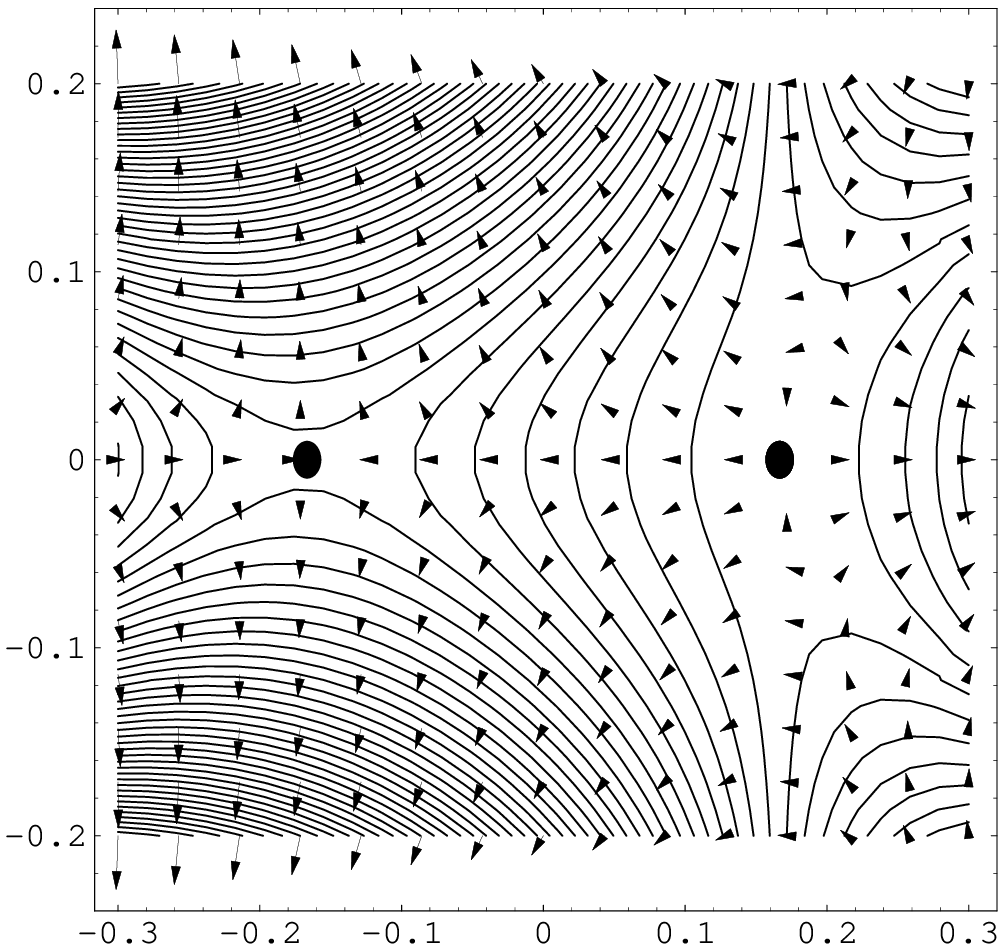}}

& 

\leavevmode\lower0in\hbox{
\includegraphics[
width=2.0211in,height=2.0211in]%
{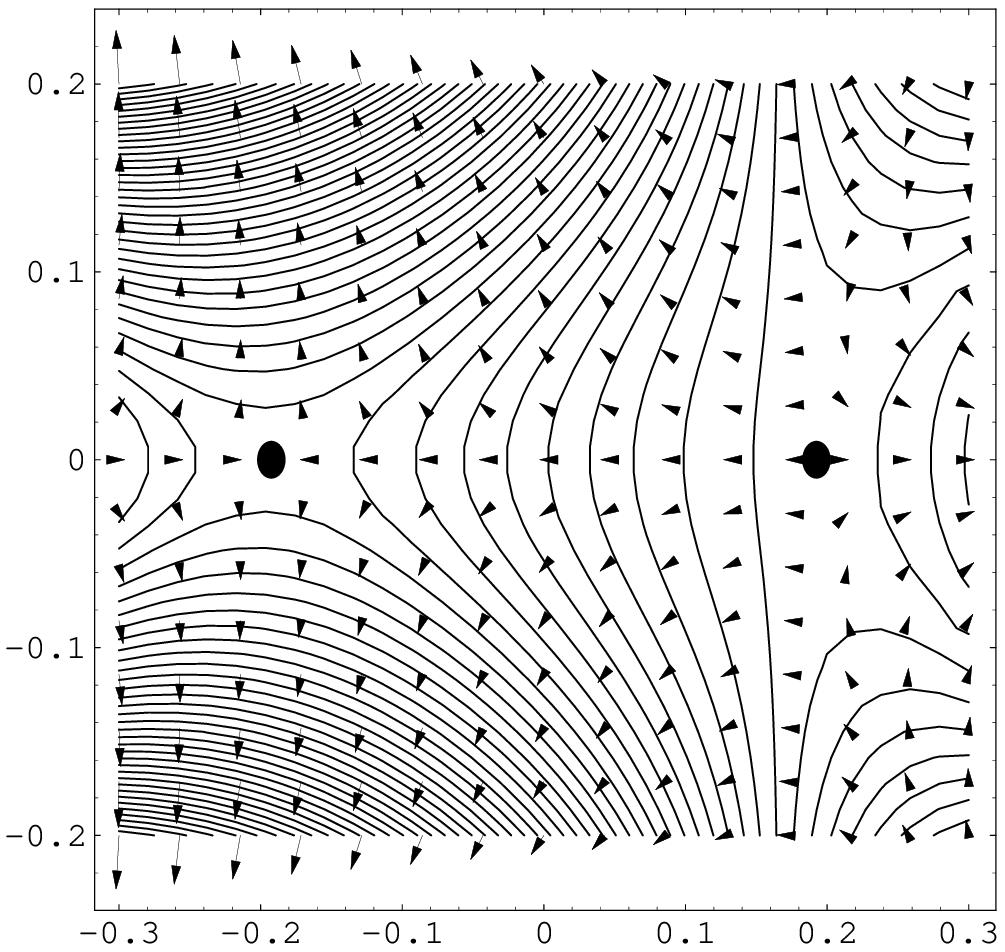}}

\\ 
\multicolumn{2}{c}{Planche 3}
\end{tabular}

\begin{figure}[tbp]
\begin{center}
\includegraphics[
width=1.4581in,height=2.7916in]%
{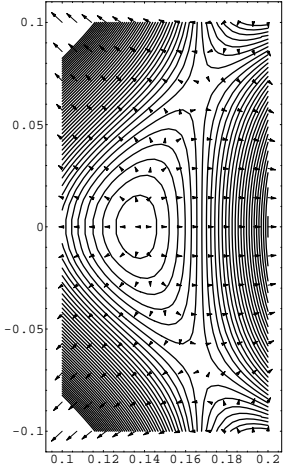}
\end{center}
\caption{Zoom sur la valeur de l'\'{e}chelle $s=(2/3)1/72$ proche de la valeur
critique $s=1/72$.}
\label{cpvfZoom}%
\end{figure}

La figure \ref{cpvfZoom} zoome sur la situation $s=(2/3)1/72$ proche de la
valeur critique $s=1/72$.

On voit bien la suite des bifurcations.\ Au d\'{e}part ($s<0$) il n'y a que
les deux cols $pc_{1\pm }$ sym\'{e}triques par rapport \`{a} l'axe des $x$
et d'abscisse fixe $1/6$. Ils se rapprochent de l'axe des $x$ comme $\sqrt{%
1-72s}$. Pour $s=0$ appara\^{i}ssent sur l'axe des $x$ un sommet $pc_{2+}$
d'abcisse $\sqrt{2s}$ et un col $pc_{2-}$ d'abscisse $-\sqrt{2s}$. Pour $%
s=1/72$, $pc_{1+}$ et $pc_{1-}$ fusionnent sur l'axe des $x$ en $x=1/6$ au
moment o\`{u} $pc_{2+}$ est pr\'{e}cis\'{e}ment en $x=\sqrt{2s}=\sqrt{2/72}%
=1/6$. La fusion de ces trois points critiques transforme $pc_{2+}$ en un
col qui s'\'{e}loigne sur l'axe des $x$ en $\sqrt{2s}$.

\begin{figure}[tbp]
\begin{center}
\includegraphics[
width=2.0211in,height=1.2496in]%
{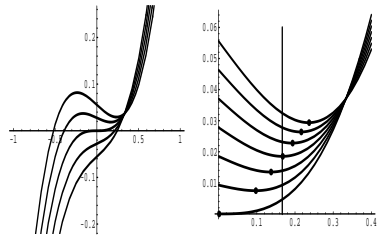}
\end{center}
\caption{Section des surfaces $f(x,y,s)$ le long de l'axe
de sym\'{e}trie $y=0$.}
\label{Medianes}%
\end{figure}

On voit clairement sur ces images l'apparition temporaire du maximum. Si le
principe de causalit\'{e} n'est pas viol\'{e} c'est \`{a} cause de la fa\c{c}%
on dont les surfaces se positionnent les unes par rapport aux autres.\ La
nappe montante reste \`{a} peu pr\`{e}s fixe alors que la nappe descendante
s'abaisse comme le montre la figure \ref{Medianes} qui repr\'{e}sente 
la section le long de l'axe
de sym\'{e}trie $y=0$.

\subsection{Plongement dans le d\'{e}ploiement de l'ombilic elliptique}

La famille $f(x,y,s)=x^{3}-6xy^{2}+y^{2}-6sx+2s$ est une d\'{e}formation de
la singularit\'{e} ombilic elliptique $g_{0}(x,y)=x^{3}-6xy^{2}$ et il est
donc pertinent de la positionner par rapport au d\'{e}ploiement universel
standard de ce dernier, $g_{w,u,v}=x^{3}-6xy^{2}+wx^{2}+ux+vy$. En fait
c'est ce qui se passe au voisinage du point critique
d\'{e}g\'{e}n\'{e}r\'{e} triple $(1/6,0)$ pour $s=1/72$ qui est le plus
int\'{e}ressant. On a

\[
f\left( \frac{1}{6}+x,y,s\right) =x^{3}-6xy^{2}+\frac{x^{2}}{2}+\left( \frac{%
1}{12}-6s\right) x+s+\frac{1}{216} 
\]

\noindent ce qui, pour $s=1/72$, donne la singularit\'{e} (centre
organisateur)

\[
f\left( \frac{1}{6}+x,y,\frac{1}{72}\right) =x^{3}-6xy^{2}+\frac{x^{2}}{2}+%
\frac{1}{54}. 
\]

Le d\'{e}ploiement universel $g_{w,u,v}=x^{3}-6xy^{2}+wx^{2}+ux+vy$ a pour
discriminant dans le $\Bbb{R}^{3}$ $(u,v,w)$ une surface conique sur une
hypocyclo\"{i}de \`{a} $3$ rebroussements. Sa section plane $w=1/2$ est une
telle hypocyclo\"{i}de $H$ d'\'{e}quation

\[
(-1+u)u^{3}+(486-648u+144u^{2})v^{2}+5184v^{4}=0 
\]

\begin{figure}[tbp]
\begin{center}
\includegraphics[
width=1.2427in,height=2.0211in]%
{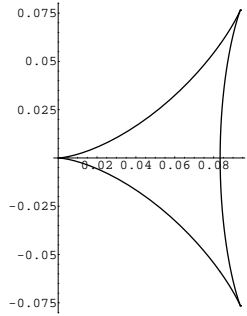}
\end{center}
\caption{La section du d\'{e}ploiement universel de l'ombilic
elliptique $g_{w,u,v}=x^{3}-6xy^{2}+wx^{2}+ux+vy$ pour $w=1/2$.}
\label{Discrim}%
\end{figure}

\noindent Elle est repr\'{e}sent\'{e}e \`{a} la figure \ref{Discrim}.

\begin{figure}[tbp]
\begin{center}
\includegraphics[
width=1.6942in,height=2.2157in]%
{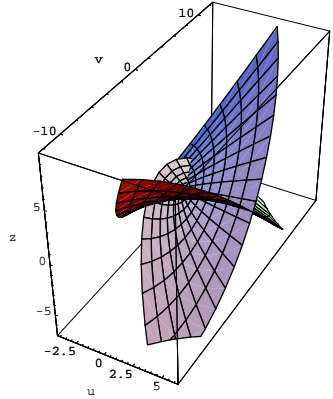}
\end{center}
\caption{Un point de vue sur le graphe des valeurs critiques de
l'ombilic elliptique $g_{1/2,u,v}(x,y)$.}
\label{GraphOv}%
\end{figure}

\begin{figure}[tbp]
\begin{center}
\includegraphics[
width=1.9865in,height=2.0211in]%
{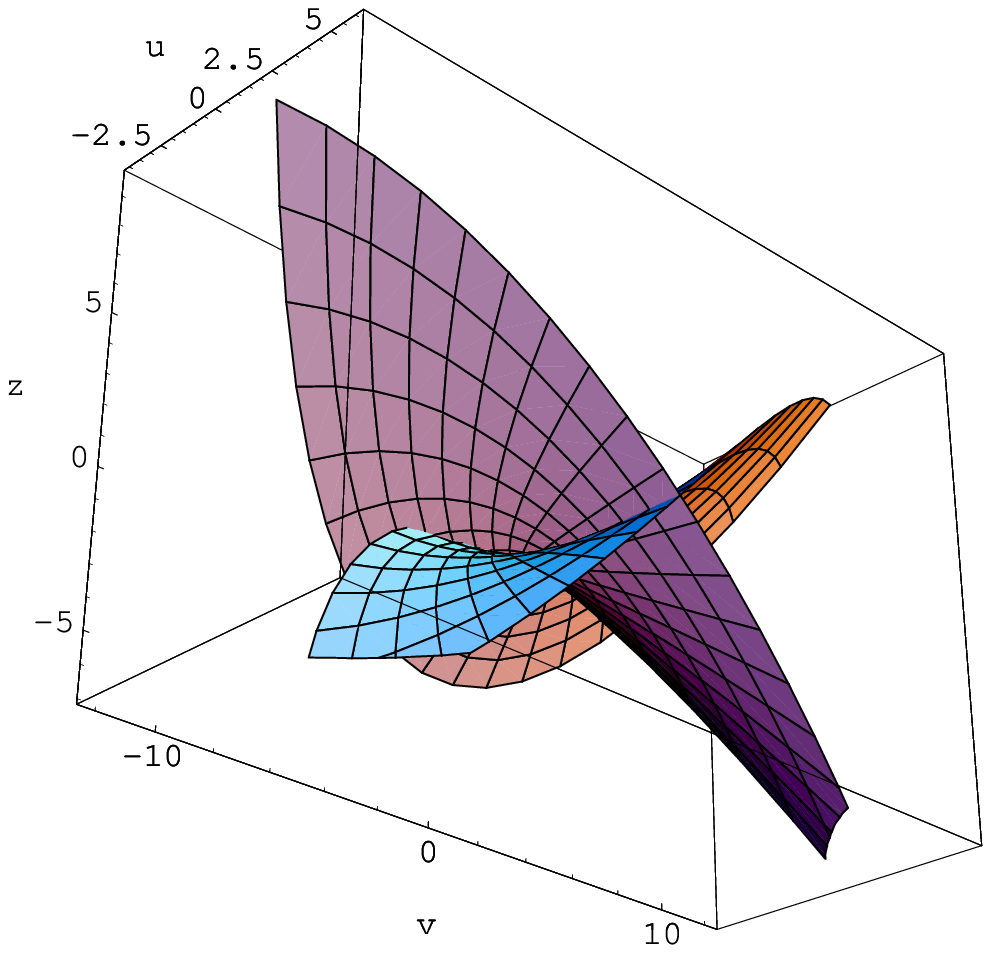}
\end{center}
\caption{Un autre point de vue sur le graphe des valeurs critiques
de l'ombilic elliptique $g_{1/2,u,v}(x,y)$.}
\label{GraphUnd}%
\end{figure}

Pour se faire une id\'{e}e de la structure de l'ombilic elliptique on peut
consid\'{e}rer dans le $\Bbb{R}^{3}$ $(u,v,z)$ le graphe des valeurs
critiques de $g_{1/2,u,v}(x,y)$. Les figures \ref{GraphOv} et \ref{GraphUnd}
en donnent deux points de vue. 

Pour visualiser le d\'{e}tail de la structure pr\`{e}s de l'origine, on a
repr\'{e}sent\'{e} \`{a} la planche $4$ trois sections au-dessus de trois
fen\^{e}tres du plan de contr\^{o}le $(u,v)$.

\begin{tabular}{ccc}

\leavevmode\lower0in\hbox{
\includegraphics[
width=1.5783in,height=1.6224in]%
{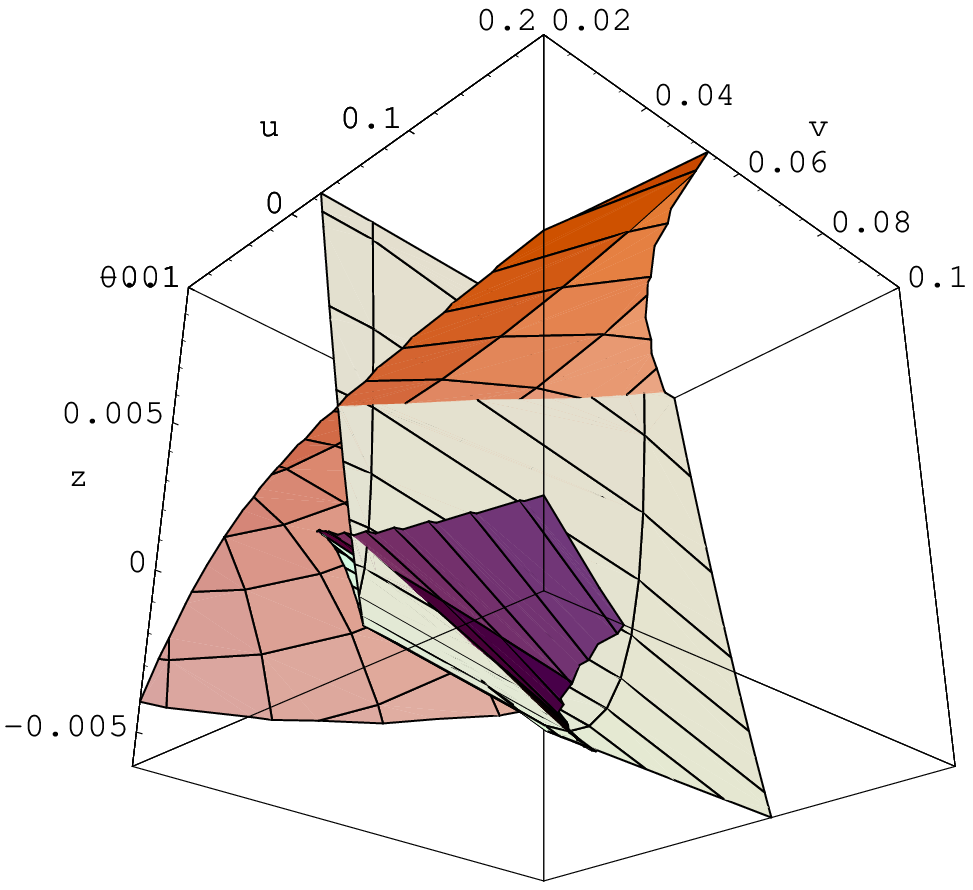}}

&

\leavevmode\lower0in\hbox{
\includegraphics[
width=1.5783in,height=1.6224in]%
{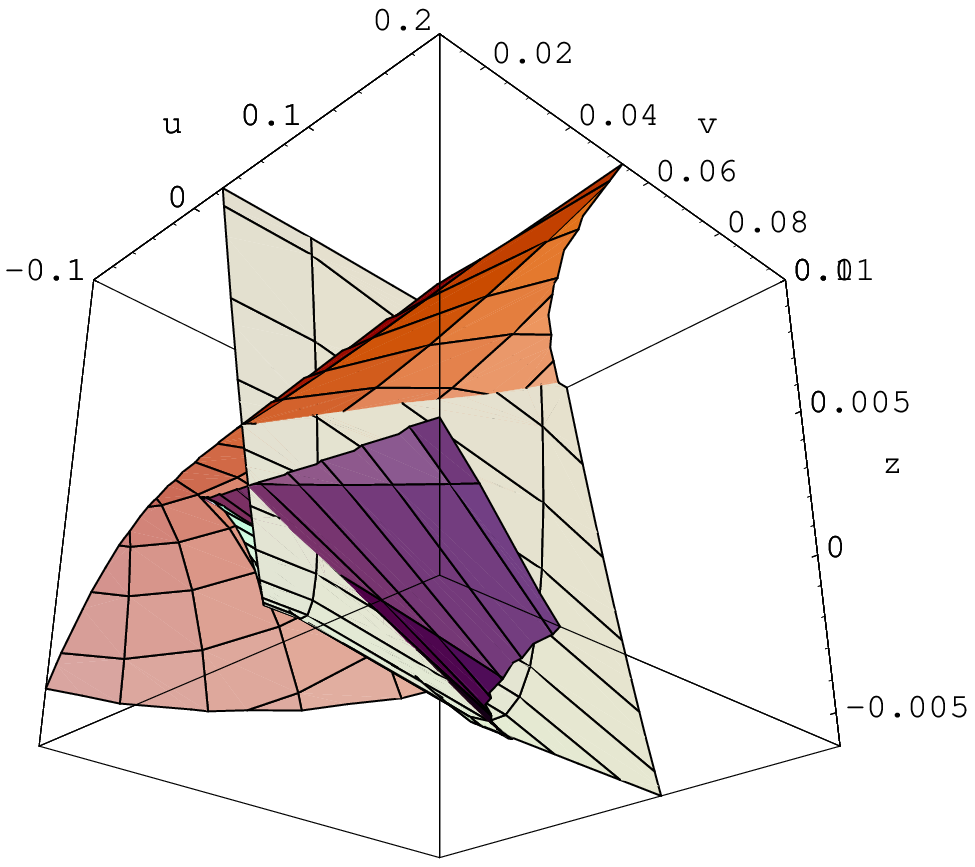}}

& 

\leavevmode\lower0in\hbox{
\includegraphics[
width=1.5783in,height=1.6224in]%
{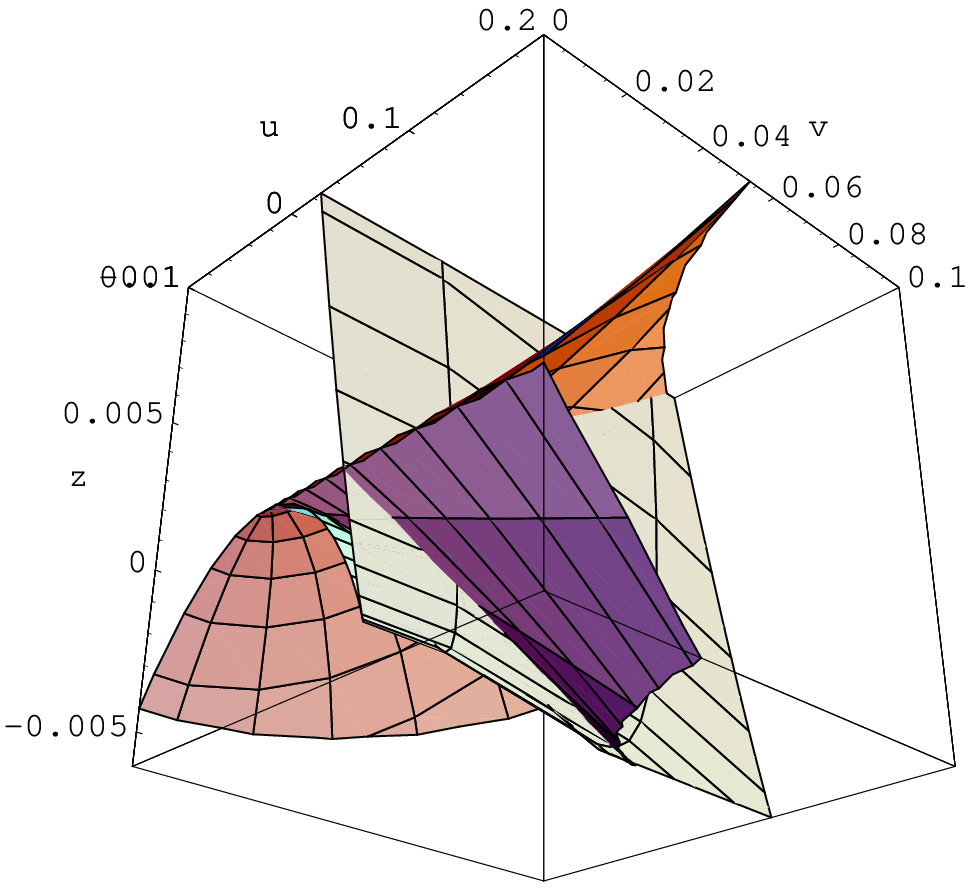}}

\\ 

\leavevmode\lower0in\hbox{
\includegraphics[
width=1.6509in,height=0.9807in]%
{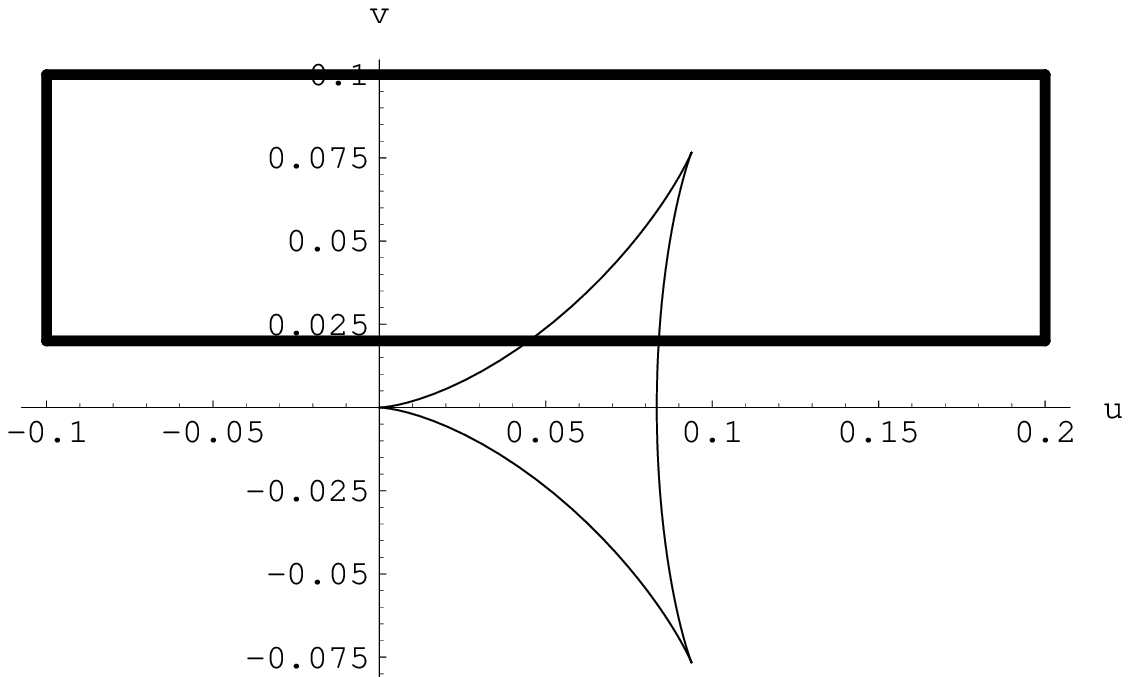}}

& 

\leavevmode\lower0in\hbox{
\includegraphics[
width=1.6509in,height=0.9807in]%
{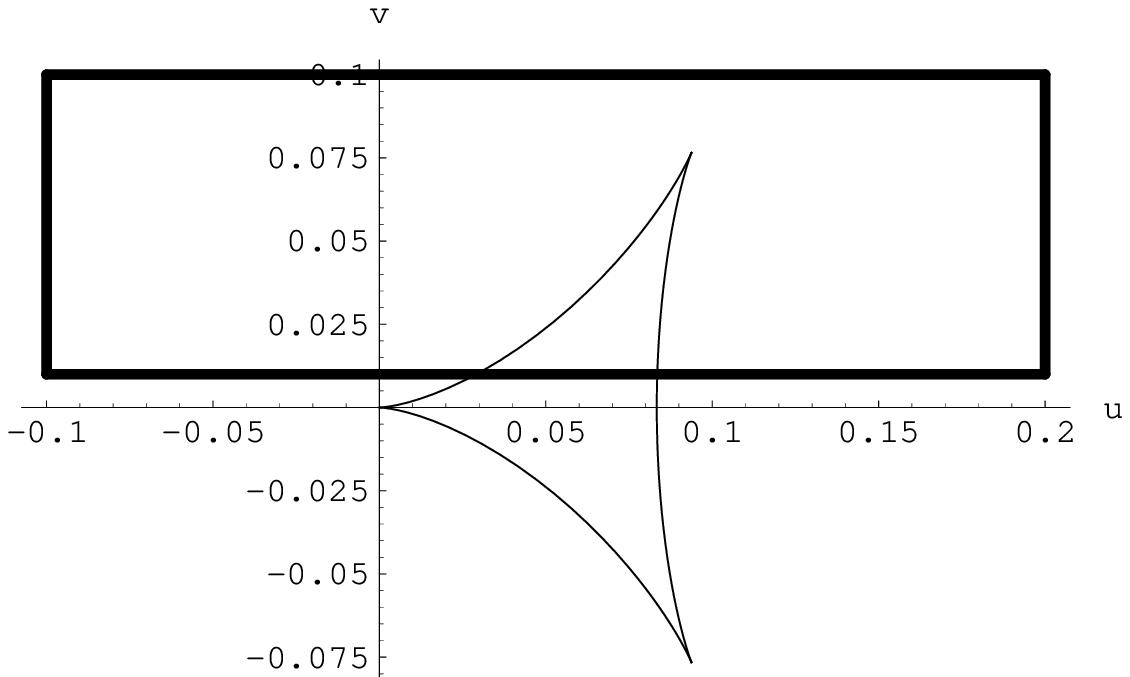}}

& 

\leavevmode\lower0in\hbox{
\includegraphics[
width=1.6509in,height=0.9807in]%
{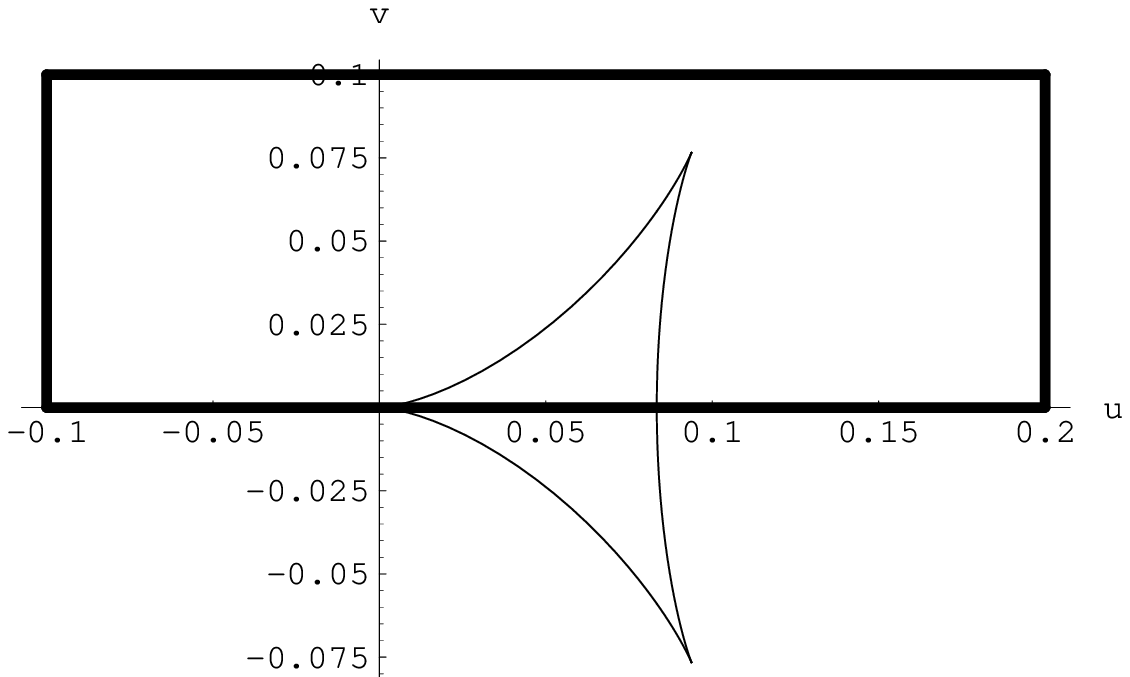}}

\\ 
\multicolumn{3}{c}{Planche $4$}
\end{tabular}

\begin{figure}[tbp]
\begin{center}
\includegraphics[
width=2.1482in,height=1.8654in]%
{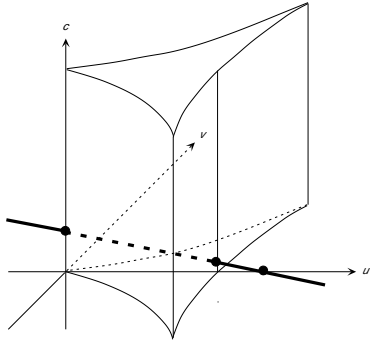}
\end{center}
\caption{Le positionnement de la solution $f(x,y,s)$ de l'\'{e}quation 
de la chaleur dans le d\'{e}ploiement de l'ombilic elliptique.}
\label{InDiscrim}%
\end{figure}

On voit que la solution \textit{stable} de l'\'{e}quation de la chaleur $%
f(x,y,s)$ correspond \`{a} un d\'{e}ploiement \`{a} un param\`{e}tre
hautement \textit{non g\'{e}n\'{e}rique} de l'ombilic elliptique. Dans le $%
\Bbb{R}^{3}$ $(u,v,c)$ le discriminant du d\'{e}ploiement $%
g_{w,u,v}=x^{3}-6xy^{2}+\frac{1}{2}x^{2}+ux+vy+c$ est le cylindre de base $H$%
. La solution stable $f(x,y,s)$ correspond \`{a} la droite $u=\left( \frac{1%
}{12}-6s\right) $, $c=s+\frac{1}{216}$ au-dessus de l'axe $u$ qui rentre
dans $H$ pour $s=0$ en $u=12$ \`{a} la hauteur $c=1/216$ et en ressort par
le cusp $u=0$, $c=1/54$ pour $s=1/72$ (cf. Figure \ref{InDiscrim}).

\section{Equations de diffusion non lin\'{e}aires et traitement d'images}

Nous avons vu comment la th\'{e}orie des singularit\'{e}s de
Whitney-Thom-Mather-Arnold pouvait \^{e}tre adapt\'{e}e aux solutions d'une
EDP parabolique de diffusion comme l'\'{e}quation de la chaleur. Mais il ne
s'agit l\`{a} que du d\'{e}but de l'histoire.\ En effet les m\'{e}thodes de
diffusion simple (homog\`{e}ne et isotrope) font probl\`{e}me car le lissage
(Gaussian blurring) qu'elles effectuent est indiff\'{e}rent \`{a} la
g\'{e}om\'{e}trie de l'image. D'o\`{u} l'id\'{e}e qui s'est progressivement
impos\'{e}e au cours des ann\'{e}es 80 que, pour pouvoir effectuer une bonne
analyse morphologique des images, il faut concilier deux exigences
apparemment contradictoires~:

\begin{description}
\item[(i)]  r\'{e}gulariser le signal de fa\c{c}on multi\'{e}chelle par
diffusion;

\item[(ii)]  pr\'{e}server la g\'{e}om\'{e}trie de l'image, c'est-\`{a}-dire
les discontinuit\'{e}s d'\'{e}l\'{e}ments diff\'{e}rentiels poss\'{e}dant
une signification g\'{e}om\'{e}trique intrins\`{e}que.
\end{description}

Pour ce faire, il faut \textit{adapter} l'\'{e}quation de diffusion \`{a} la
pr\'{e}servation de ces \'{e}l\'{e}ments diff\'{e}rentiels construits \`{a}
partir des jets successifs de la fonction $I_{s}(x,y)$.

Le cas le plus simple est celui des \textit{bords} qui d\'{e}limitent les
domaines homog\`{e}nes d'une image et sont essentiels \`{a} la
d\'{e}finition de ses constituants. Un bord est id\'{e}alement une
discontinuit\'{e} du gradient $\nabla I_{s}(x,y)$ de $I_{s}$. Pour qu'une
diffusion pr\'{e}serve le caract\`{e}re discontinu des bords tout en les
simplifiant progressivement, il faut qu'elle soit \textit{anisotrope} et,
plus pr\'{e}cis\'{e}ment, inhib\'{e}e dans la direction du gradient. L'EDP
la plus simple poss\'{e}dant cette propri\'{e}t\'{e} est l'\'{e}quation de
diffusion non lin\'{e}aire~: 
\[
\frac{\partial I_{s}}{\partial s}=\left| \nabla I_{s}\right| \text{div}%
\left( \frac{\nabla I_{s}}{\left| \nabla I_{s}\right| }\right) =\Delta I_{s}-%
\frac{H\left( \nabla I_{s},\nabla I_{s}\right) }{\left| \nabla I_{s}\right|
^{2}}
\]

\noindent o\`{u} $H$ est le Hessien de $I_{s}$. Elle est uniform\'{e}ment
parabolique le long des courbes de niveau de $I_{s}$ mais totalement
d\'{e}g\'{e}n\'{e}r\'{e}e dans la direction du gradient. Elle fait
\'{e}voluer les lignes de niveau~--- et donc en particulier les bords~---
comme des fronts avec une vitesse normale \'{e}gale \`{a} leur courbure.

Jean-Michel Morel et ses coll\`{e}gues (Alvarez \textit{et al}. 1992) ont
introduit un lissage gaussien suppl\'{e}mentaire dans ces EDP de fa\c{c}on
\`{a} contr\^{o}ler la vitesse de diffusion en la couplant \`{a} la
g\'{e}om\'{e}trie de l'image. L'\'{e}quation de base devient alors~:

\[
\frac{\partial I_{s}}{\partial s}=g\left( \left| G*\nabla I_{s}\right|
\right) \left| \nabla I_{s}\right| \text{div}\left( \frac{\nabla I_{s}}{%
\left| \nabla I_{s}\right| }\right) 
\]

\noindent o\`{u} $G$ est un noyau gaussien et $g(x)$ une fonction
d\'{e}croissante telle que 
$g(x)\longrightarrow 0$ quand $x\rightarrow \infty$. 
Ils ont \'{e}galement montr\'{e} que si l'on impose une
contrainte d'invariance affine (au lieu de l'invariance euclidienne), l'EDP
type devient~:

\begin{center}
\[
\frac{\partial I_{s}}{\partial s}=\left| \nabla I_{s}\right| \text{div}%
\left( \frac{\nabla I_{s}}{\left| \nabla I_{s}\right| }\right) ^{1/3}. 
\]
\end{center}

On peut g\'{e}n\'{e}raliser ces \'{e}quations de fa\c{c}on \`{a}
am\'{e}liorer l'analyse morphologi\-que de l'image. Par exemple les
``cr\^{e}tes'' d'une forme~--- qui fournissent une ``squel\'{e}tisation''
essentielle \`{a} son analyse morphologique, ce que Ren\'{e} Thom avait
souvent soulign\'{e} \`{a} la suite de Harry Blum en utilisant le concept de
``cut-locus'')~--- peuvent \^{e}tre extraites du signal en tant que des
discontinuit\'{e}s de la direction du gradient qui sont pr\'{e}serv\'{e}es
de la diffusion.

En consid\'{e}rant les lignes de niveau des fonctions $I_{s}(x,y)$ on
obtient des \'{e}volutions de courbes planes $C_{s}$ ferm\'{e}es qui se
propagent comme des fronts conform\'{e}ment \`{a} une loi du type $\frac{%
\partial p}{\partial s}=F(\kappa )\vec{n}$ o\`{u} $p$ est un point de $C_{s}$%
, $\vec{n}$ la normale (externe) en $p$ \`{a} $C_{s}$ et $\kappa $ la
courbure de $C_{s}$ en $p$ . Les cas les plus \'{e}tudi\'{e}s sont~:

\begin{description}
\item[(i)]  celui de la propagation \`{a} vitesse constante~: $\frac{%
\partial p}{\partial s}=\pm \vec{n},\frac{\partial \kappa }{\partial s}%
=-\kappa ^{2}$ (mod\`{e}les de propagation d'ondes de type ``grassfire''
engendrant le cut-locus de la courbe), et

\item[(ii)]  celui de la propagation \`{a} une vitesse proportionnelle \`{a}
la courbure ($\tau $ est l'abscisse curviligne de la courbe)~: $\frac{%
\partial p}{\partial s}=-\kappa \vec{n},\frac{\partial \kappa }{\partial s}=%
\frac{\partial ^{2}\kappa }{\partial \tau ^{2}}+\kappa ^{3}$

\item[(iii)]  Stanley Osher et James Sethian (1988) ont \'{e}galement
\'{e}tudi\'{e} les cas interm\'{e}diaires o\`{u} $F(\kappa )=1-\varepsilon
\kappa $. La courbure $\kappa $\textit{\ }satisfait alors une \'{e}quation
de r\'{e}action-diffusion du type~: $\frac{\partial \kappa }{\partial s}%
=\varepsilon \frac{\partial ^{2}\kappa }{\partial \tau ^{2}}+\varepsilon
\kappa ^{3}-\kappa ^{2}.$
\end{description}

Sous les titres de ``curve shortening'', ``flow by curvature'' ou ``heat
flow on isometric immersions'', le cas (ii) a \'{e}t\'{e}
particuli\`{e}rement investigu\'{e} par des g\'{e}om\`{e}tres comme M. Gage,
R. Hamilton (1986), M. Grayson (1987), S. Osher, J.\ Sethian (1988), J.
Sethian (1990), L.C. Evans et J. Spruck (1991).\footnote{%
En fait, la th\'{e}orie vient de Richard Hamilton qui cherchait \`{a}
r\'{e}soudre des probl\`{e}mes de Relativit\'{e} g\'{e}n\'{e}rale. En
utilisant l'\'{e}quation de la chaleur, il a montr\'{e} que si $X$ est une
vari\'{e}t\'{e} riemannienne compacte de dimension $3$ avec courbure de
Ricci $R_{ij}>0$, alors $X$ admet une m\'{e}trique riemannienne \`{a}
courbure $>0$ constante. Or, ces derni\`{e}res sont classifi\'{e}es. Il
cherchait \'{e}galement \`{a} engendrer des g\'{e}od\'{e}siques ferm\'{e}es
\`{a} partir de courbes ferm\'{e}es quelconques. (Cf. le s\'{e}minaire
Bourbaki Bourguignon [1985]).} Si $j_{s}:S^{1}\rightarrow \Bbb{R}^{2}$ est
l'immersion isom\'{e}trique d\'{e}finissant $C_{s}$, comme on a $\Delta %
j_{s}=-\kappa \vec{n}$, l'\'{e}quation de diffusion $\frac{\partial j_{s}}{%
\partial s}=-\kappa \vec{n}$ est en fait l'\'{e}quation de la chaleur (pour
les immersions) $\frac{\partial j_{s}}{\partial s}=\Delta j_{s}$. Dans
l'espace fonctionnel $\mathcal{J}$ des immersions $j:S^{1}\rightarrow \Bbb{R}%
^{2}$, cette \'{e}quation d\'{e}finit le champ de gradient de la fonction de
Morse donnant la longueur de la courbe image $C=j(S^{1})$.

Un th\'{e}or\`{e}me fondamental de ``curve shortening'' d\^{u} \`{a} M. Gage
et R. Hamilton (1985) et g\'{e}n\'{e}ralis\'{e} par Matthew Grayson (1987)
dit que si $C_{0}$ est une courbe plong\'{e}e dans le plan (m\^{e}me
tr\`{e}s sinueuse), alors l'\'{e}quation de la chaleur la contracte sur un
point, $C_{s}$ devenant asymptotiquement un cercle (convergence pour la
norme $C^{\infty }$). En particulier $C_{s}$ devient convexe avant de
pouvoir d\'{e}velopper des singularit\'{e}s.

\section{Conclusion}

C'est en fait toute la conception morphodynamique thomienne sur le r\^{o}le
des singularit\'{e}s, des bifurcations et des ruptures de sym\'{e}tries dans
la perception qui se trouve confirm\'{e}e par un nombre croissant de travaux
contemporains en vision computationnelle et en neurosciences de la vision.
Citons pour conclure trois exemples~:

\begin{description}
\item[(i)]  Celui des processus de segmentation des images en r\'{e}gion
homog\`{e}nes d\'{e}limit\'{e}es par des discontinuit\'{e}s qualitatives.\
Ce probl\`{e}me central \'{e}tait \`{a} la base des mod\`{e}les
morphologiques de Thom d\`{e}s le d\'{e}but avec la diff\'{e}rence entre les
points r\'{e}guliers d'un substrat o\`{u} toutes les qualit\'{e}s
ph\'{e}nom\'{e}nologiques (couleur, texture, etc.) sont localement
homog\`{e}nes et les points singuliers o\`{u} au moins une qualit\'{e}
ph\'{e}nom\'{e}nologique subit une discontinuit\'{e} qualitative. Il a
\'{e}t\'{e} formul\'{e} comme un (difficile) probl\`{e}me variationnel de
type ``free boundary problem'' par David Mumford dans les ann\'{e}es 80 et
\`{a} donn\'{e} lieu \`{a} tout un ensemble de travaux, en particulier
autour d'Ennio De Giorgi \`{a} Pise (Luigi Ambrosio, Gianni Dal Maso). La
conjecture de Mumford sur la nature des points singuliers des solutions
optimales a presque \'{e}t\'{e} d\'{e}montr\'{e}e par Alexis Bonnet et Guy
David.

\item[(ii)]  La squel\'{e}tisation des formes planes au moyen de leur cut
locus, c'est-\`{a}-dire des singularit\'{e}s d'un processus de propagation
\`{a} la Huyghens \`{a} partir du bord.\ Cette id\'{e}e profond\'{e}ment
morphologique a \'{e}t\'{e} introduite dans les ann\'{e}es 60 par le
sp\'{e}cialiste de la vision Harry Blum [1973] et a toujours \'{e}t\'{e}
tr\`{e}s fortement d\'{e}fendue par Ren\'{e} Thom. L'avantage du squelette
d'une forme est d'\^{e}tre un graphe de dimension 1 constitu\'{e} d'arcs se
rejoignant g\'{e}n\'{e}riquement en des points triples ou s'arr\^{e}tant en
des points singuliers. Il fournit donc automatiquement une d\'{e}composition
qualitative de la forme en constituants qui sont des cylindres
g\'{e}n\'{e}ralis\'{e}s (cf.\ Marr [1982]). L'hom\'{e}omorphismes des
squelettes d\'{e}finit une relation d'\'{e}quivalence qui est plus forte que
les diff\'{e}omorphismes et plus faible que les isom\'{e}tries tout en
\'{e}tant tr\`{e}s diff\'{e}rente de l'\'{e}quivalence conforme. Pendant
longtemps cette m\'{e}thode \'{e}tait peu pris\'{e}e car instable et trop
sensible aux d\'{e}formations du bord, toute petite asp\'{e}rit\'{e} du bord
engendrant une branche suppl\'{e}mentaire du squelette.\ Mais le
d\'{e}veloppement de m\'{e}thodes de pruning multi\'{e}chelle a rendu la
m\'{e}thode stable et elle est aujourd'hui tr\`{e}s utilis\'{e}e. Un des
meilleurs sp\'{e}cialistes en est Benjamin Kimia.

\item[(iii)]  La reconstruction des surfaces plong\'{e}es dans $\Bbb{R}^{3}$
\`{a} partir de la famille de leurs contours apparents.\ L\`{a} encore
Ren\'{e} Thom a toujours insist\'{e} sur l'importance du probl\`{e}me. Si $S$
est une surface dans $\Bbb{R}^{3}$, \`{a} chaque point de vue $p$ (plan de
projection et direction de projection, les $p$ forment une grassmannienne $%
\mathcal{P}$) est associ\'{e} un contour apparent $C_{p}$.\ Le type
qualitatif de $C_{p}$ stratifie $\mathcal{P}$ \`{a} travers une hypersurface 
$H$ et le probl\`{e}me est de comprendre l'\'{e}quivalence entre la
g\'{e}om\'{e}trie de $S$ et la g\'{e}om\'{e}trie de la stratification
d\'{e}finie par $H$. Par exemple si l'on conna\^{i}t des \'{e}quations
alg\'{e}briques pour $S$ peut-on en d\'{e}duire des \'{e}quations pour $H$?\
Le probl\`{e}me est difficile. Or le syst\`{e}me visuel est capable
d'anticiper avec une pr\'{e}cision remarquable les variations de $C_{p}$
lorsque l'on tourne autour de $S$. On ne conna\^{i}t pas les algorithmes
neuralement impl\'{e}ment\'{e}s sous-jacents \`{a} cette performance
perceptive.\footnote{%
Le lecteur int\'{e}ress\'{e} par l'actualit\'{e} de certains mod\`{e}les de
g\'{e}om\'{e}trie diff\'{e}rentielle en vision pourra consulter le double
num\'{e}ro sp\'{e}cial (97, 2-3) du \textit{Journal of Physiology-Paris}
(Petitot, Lorenceau [2003], disponible sur www.sciencedirect.com).}
\end{description}

Tous ces travaux convergents montrent que le r\^{e}ve de Ren\'{e} Thom de
d\'{e}velopper une ``physique ph\'{e}nom\'{e}nologique'' et une
``s\'{e}miophysique'' de la perception est en train de se r\'{e}aliser. Il
en va de m\^{e}me dans le domaine de la morphogen\`{e}se biologique o\`{u}
Ren\'{e} Thom estimait avoir introduit ses id\'{e}es les plus fondamentales.
L\`{a} aussi, comme pour les sciences cognitives, on peut dire que
l'utilisation de mod\`{e}les dynamiques et de la th\'{e}orie des
bifurcations engendrant des discontinuit\'{e}s qualitatives dans des
substrats physiques, s'est beaucoup d\'{e}velopp\'{e}e parall\`{e}lement
\`{a} d'autres mod\`{e}les de th\'{e}orie des patterns comme les mod\`{e}les
de r\'{e}action-diffusion de type Turing. Ce qui a peut-\^{e}tre le plus
chang\'{e} d'avec la fa\c{c}on dont Ren\'{e} Thom lui-m\^{e}me voyait les
choses, est qu'il semble d\'{e}sormais possible d'\'{e}tablir une
convergence avec les r\'{e}sultats de la g\'{e}nomique. L'on commence \`{a}
comprendre le contr\^{o}le g\'{e}n\'{e}tique du d\'{e}veloppement et des
r\'{e}actions morphog\'{e}n\'{e}tiques qui int\'{e}ressaient Ren\'{e} Thom.

Ren\'{e} Thom \'{e}tait non seulement un puissant g\'{e}nie math\'{e}matique
et un ma\^{i}tre humainement extraordinaire, mais aussi l'un des esprits
philosophiques les plus inspir\'{e}s du XXe si\`{e}cle. Ce f\^{u}t un
privil\`{e}ge inoubliable que d'avoir pu collaborer \'{e}troitement et
contin\^{u}ment avec lui. Au-del\`{a} des al\'{e}as humains trop humains de
la sociologie des milieux acad\'{e}miques, ses id\'{e}es ont d'ores et
d\'{e}j\`{a} triomph\'{e} dans de nombreux domaines et, j'en suis convaincu,
continueront \`{a} le faire, ainsi que toutes les grandes id\'{e}es
appartenant au patrimoine de la pens\'{e}e.

\section{Bibliographie}
\ 
\indent ALVAREZ, L., LIONS, P.L., MOREL, J.M., 1992. ``Image selective smoothing and
edge detection by non linear diffusion, \textit{SIAM J. Numer. Anal.}, 29,
845-866.

BLUM, H., 1973. ``Biological Shape and Visual Science'', \textit{Journal of
Theoretical Biology}, 38, 205-287.

BOURGUIGNON, J.-P., 1985. ``L'\'{e}quation de la chaleur associ\'{e}e \`{a}
la courbure de Ricci'', \textit{S\'{e}minaire Bourbaki 653}, Ast\'{e}risque.

CHENCINER, A., 1973. ``Travaux de Thom et Mather sur la Stabilit\'{e}
Topologique'', \textit{S\'{e}minaire Bourbaki 424}.

CHENCINER, A., 1980. ``Singularit\'{e}s des Fonctions Diff\'{e}rentiables'', 
\textit{Encyclop\ae dia Universalis}, Paris.

DAMON, J., 1995. ``Local Morse Theory for Solutions to the Heat Equation and
Gaussian Blurring'', \textit{Journal of Differential Equations}, 115, 2,
368-401.

EVANS, L.C., SPRUCK, J., 1991. ``Motion of Level Sets by Mean Curvature.
I'', \textit{J. Differential Geometry}, 33, 635-681.

FLORACK, L., 1993. \textit{The Syntactical Structure of Scalar Images},
Th\`{e}se, Universit\'{e} d'Utrecht.

GAGE, M., HAMILTON, R. S., 1986. ``The Heat equation shrinking convex plane
curves'', \textit{J. Differential Geometry}, 25, 69-96.

GIBSON, J.J., 1979. \textit{The Ecological Approach to Visual Perception},
Boston, Houghton-Mifflin.

GOLUBITSKY, M., GUILLEMIN, V., 1973. \textit{Stable Mappings and their
Singularities}, Graduate Texts in Mathematics, 14, Springer, New-York,
Heidelberg, Berlin.

GRAYSON, M., 1987. ``The Heat Equation Shrinks Embedded Plane Curves to
Round Points'', \textit{J. Differential Geometry}, 26, 285-314.

KOENDERINK, J.J., 1984. ``The Stucture of Images'', \textit{Biological
Cybernetics}, 50, 363-370.

KOENDERINK, J.J., VAN DOORN, A.J., 1986. ``Dynamic Shape'', \textit{%
Biological Cybernetics}, 53, 383-396.

MARR, D., 1982. \textit{Vision}, San Francisco, Freeman.

MATHER, J., 1968(a). ``Stability of $C^{\infty }$ Mappings I~: the Division
Theorem``, \textit{Ann. of Math}., 87, 89-104.

MATHER, J., 1968(b). ``Stability of $C^{\infty }$ Mappings III~: Finitely
Determined Map Germs'', \textit{Publications Math\'{e}matiques de l'IHES},
35, 127-156, Presses Universitaires de France, Paris.

MATHER, J., 1969(a). ``Stability of $C^{\infty }$ Mappings II~:
Infinitesimal Stability Implies Stability'', \textit{Ann. of Math}., 89,
254-291.

MATHER, J., 1969(b). ``Stability of $C^{\infty }$ Mappings IV~:
Classification of Stable Germs by R -algebras'', \textit{Publications
Math\'{e}matiques de l'IHES}, 37, 223-248.

MATHER, J., 1970. ``Stability of $C^{\infty }$ Mappings V~:
Transversality'', \textit{Advances in Math}., 4, 301-336.

MATHER, J., 1971. ``Stability of $C^{\infty }$ Mappings VI~: the Nice
Dimensions'', \textit{Proceedings of Liverpool Singularities Symposium}
(C.T.C. Wall, ed.), Lecture Notes in Mathematics 192, Springer, New-York,
Heidelberg, Berlin., 207-253.

MUMFORD, D., SHAH, J., 1988. ``Boundary Detection by Minimizing
Functionals'', \textit{Proceedings IEEE Computer Vision and Pattern
Recognition Conference}, Ann Arbor, Michigan.

OSHER, S., SETHIAN, J. A., 1988. ``Fronts Propagating with
Curvature-Dependent Speed~: Algorithms Based on Hamilton-Jacobi
Formulations'', \textit{Journal of Computational Physics}, 79, 1, 12-49.

PETITOT, J., 1989. ``Forme'', \textit{Encyclop\ae dia Universalis}, XI,
712-728, Paris.

PETITOT, 1989 (ed.). \textit{Logos et Th\'{e}orie des Catastrophes}
(Colloque de Cerisy \`{a} partir de l'oeuvre de Ren\'{e} Thom), Gen\`{e}ve,
Editions Patino.

PETITOT, J., 1992. \textit{Physique du Sens}, Paris, Editions du CNRS.

PETITOT, J., 1994. ``La s\'{e}miophysique~: de la physique qualitative aux
sciences cognitives'', \textit{Passion des Formes, \`{a} Ren\'{e} Thom} (M.
Porte \'{e}d.), 499-545, E.N.S. Editions Fontenay-Saint Cloud.

PETITOT, J., 2003 (ed. avec Jean Lorenceau). \textit{Neurogeometry and
Visual Perception}, \textit{Journal of Physiology-Paris}, 97, 2-3.

SETHIAN, J.A., 1990. ``Numerical Algorithms for Propagating Interfaces~:
Hamilton-Jacobi Equations and Conservation Laws'', \textit{J. Differential
Geometry}, 31, 131-161.

THOM, R., 1972. \textit{Stabilit\'{e} structurelle et Morphogen\`{e}se}, New
York, Benjamin, Paris, Ediscience.

THOM, R., 1980. \textit{Mod\`{e}les math\'{e}matiques de la Morphogen\`{e}se}
(2\`{e}me ed.), Paris, Christian Bourgois.

WITKIN, A., 1983. ``Scale-Space Filtering'', \textit{Proc. Int. Joint Conf.
on Artificial Intelligence}, 1019-1021, Karlsruhe.

\end{document}